\documentclass[amsbook,12pt,a4paper]{article}
\usepackage{epsfig, amsfonts, amssymb, amsmath}
\usepackage[usenames]{color}
\usepackage{picinpar}
\usepackage{pgf,pgfarrows,pgfnodes,pgfautomata,pgfheaps,pgfshade}
\newcommand{\nl}{\mbox{}\\}

\begin{document}
\thispagestyle{empty}
%
% ----------------------------------------------------------------
%
\mbox{} \vspace{-1.000cm} \\
\begin{center}
{\Large \bf
General asymptotic supnorm estimates for} \\
\mbox{} \vspace{-0.350cm} \\
{\Large \bf
solutions of one-dimensional advection-diffusion} \\
\mbox{} \vspace{-0.350cm} \\
{\Large \bf
equations in heterogeneous media, I}
\nl
\mbox{} \vspace{-0.150cm} \\
\nl
{\bf \sc Jos\'e A. Barrionuevo,
Lucas S. Oliveira
and
Paulo R. Zingano} \\
\mbox{} \vspace{-0.350cm} \\
{\small
Departamento de Matem\'atica Pura e Aplicada} \\
\mbox{} \vspace{-0.620cm} \\
{\small
Universidade Federal do Rio Grande do Sul} \\
\mbox{} \vspace{-0.620cm} \\
{\small
Porto Alegre, RS 91509-900, Brazil} \\
\nl
\mbox{} \vspace{-0.350cm} \\
{\bf Abstract} \\
\mbox{} \vspace{-0.150cm} \\
\begin{minipage}{12.50cm}
{\footnotesize
We derive general bounds for the large time size
of supnorm values
$ {\displaystyle
\|\, u(\cdot,t) \,\|_{\mbox{}_{\scriptstyle L^{\infty}(\mathbb{R})}}
\!
} $ \linebreak
\mbox{} \vspace{-0.500cm} \\
of solutions to one-dimensional advection-diffusion equations \\
\mbox{} \vspace{-0.650cm} \\
\begin{equation}
\notag
u_t \;\!+\, (\;\! b(x,t) \;\!u \;\!)_{x} \;\!=\;
u_{xx},
\qquad x \in \mathbb{R}, \; t > 0
\end{equation}
\mbox{} \vspace{-0.530cm} \\
with initial data
$ {\displaystyle
u(\cdot,0) \in L^{p_{\mbox{}_{\!\;\!0}}}(\mathbb{R}) \cap L^{\infty}(\mathbb{R})
} $
for some $ 1 \leq p_{\mbox{}_{\!\;\!0}} \!< \infty $,
and arbitrary \linebreak
\mbox{} \vspace{-0.700cm} \\
bounded
advection speeds $ b(x,t) $,
introducing new techniques based on suitable \linebreak
\mbox{} \vspace{-0.700cm} \\
energy arguments.
Some open problems and related results
are also given. \\
}
\end{minipage}
\end{center}
\nl
\mbox{} \vspace{-0.650cm} \\
\nl
\nl
\mbox{} \vspace{-0.450cm} 
\mbox{} \hfill
heterogeneous media,
forced advection,
supnorm estimates,
large time behavior. \\
\setcounter{page}{0}
\newpage
%
%
% ----------------------------------------------------------------
%
%
\mbox{} \vspace{-2.000cm} \\
\begin{center}
{\Large \bf
General asymptotic supnorm estimates for} \\
\mbox{} \vspace{-0.350cm} \\
{\Large \bf
solutions of one-dimensional advection-diffusion} \\
\mbox{} \vspace{-0.350cm} \\
{\Large \bf
equations in heterogeneous media, I}
\nl
\mbox{} \vspace{-0.150cm} \\
\nl
{\bf \sc Jos\'e A. Barrionuevo,
Lucas S. Oliveira
and
Paulo R. Zingano} \\
\mbox{} \vspace{-0.350cm} \\
{\small
Departamento de Matem\'atica Pura e Aplicada} \\
\mbox{} \vspace{-0.620cm} \\
{\small
Universidade Federal do Rio Grande do Sul} \\
\mbox{} \vspace{-0.620cm} \\
{\small
Porto Alegre, RS 91509-900, Brazil} \\
\nl
\mbox{} \vspace{-0.350cm} \\
{\bf Abstract} \\
\mbox{} \vspace{-0.150cm} \\
\begin{minipage}{12.50cm}
{\footnotesize
We derive general bounds for the large time size
of supnorm values
$ {\displaystyle
\|\, u(\cdot,t) \,\|_{\mbox{}_{\scriptstyle L^{\infty}(\mathbb{R})}}
\!
} $ \linebreak
\mbox{} \vspace{-0.500cm} \\
of solutions to one-dimensional advection-diffusion equations \\
\mbox{} \vspace{-0.650cm} \\
\begin{equation}
\notag
u_t \;\!+\, (\;\! b(x,t) \;\!u \;\!)_{x} \;\!=\;
u_{xx},
\qquad x \in \mathbb{R}, \; t > 0
\end{equation}
\mbox{} \vspace{-0.530cm} \\
with initial data
$ {\displaystyle
u(\cdot,0) \in L^{p_{\mbox{}_{\!\;\!0}}}(\mathbb{R}) \cap L^{\infty}(\mathbb{R})
} $
for some $ 1 \leq p_{\mbox{}_{\!\;\!0}} \!< \infty $,
and arbitrary \linebreak
\mbox{} \vspace{-0.700cm} \\
bounded
advection speeds $ b(x,t) $,
introducing new techniques based on suitable \linebreak
\mbox{} \vspace{-0.700cm} \\
energy arguments.
Some open problems and related results
are also given. \\
}
\end{minipage}
\end{center}
\nl
%
% ************************************************************
% *                                                          *
% *              Section 1: Introduction                     *
% *                                                          *
% ************************************************************
%
\par
{\bf \S 1. Introduction} \\
\par
In this work,
we obtain very general large time estimates
for supnorm values of solutions $ u(\cdot,t) $
to parabolic initial value problems of the form \\
\mbox{} \vspace{-0.600cm} \\
\begin{equation}
\tag{1.1$a$}
u_t \;\!+\, (\;\!b(x,t) \;\!u \;\!)_{x}
\;\!=\;
u_{xx},
\qquad
x \in \mathbb{R}, \; t > 0,
\end{equation}
\mbox{} \vspace{-0.750cm} \\
\begin{equation}
\tag{1.1$b$}
u(\cdot,0) \,=\, u_{\mbox{}_{\!\;\!0}} \in
L^{p_{\mbox{}_{\!\;\!0}}}(\mathbb{R}) \cap L^{\infty}(\mathbb{R}),
\qquad
1 \leq p_{\mbox{}_{\!\;\!0}} \!< \infty,
\end{equation}
\mbox{} \vspace{-0.150cm} \\
for arbitrary continuously differentiable advection fields
$ \;\!b \in L^{\infty}(\mathbb{R} \times [\;\!0, \infty\;\![\;\!) $.
Here, \linebreak
by {\em solution\/} to (1.1) in some time interval
$ [\;\!0, T_{\mbox{}_{\scriptstyle \!\ast}}[ $,
$ \:\!0 < T_{\mbox{}_{\scriptstyle \!\ast}} \!\leq \infty $,
we mean a function \linebreak
$ {\displaystyle
u \!\:\!:\;\! \mathbb{R} \times [\;\!0, T_{\mbox{}_{\scriptstyle \!\ast}} [
\;\rightarrow \mathbb{R}
} $
which is bounded in each strip
$ S_{{\scriptstyle T}} \!\:\!=\, \mathbb{R} \times [\;\!0, T\:\!] $,
$ 0 < T \!\,\!< T_{\mbox{}_{\scriptstyle \!\ast}} $, \linebreak
% (i.e.,
% $ {\displaystyle
% u(\cdot,t) \in
% L^{\infty}_{\tt loc}
% (\;\![\;\!0, T_{\mbox{}_{\scriptstyle \!\ast}} [, L^{\infty}(\mathbb{R}))
% } $),
solves equation (1.1$a$) in the classical sense for
$ \;\!0 < t < T_{\mbox{}_{\scriptstyle \!\ast}} $,
and satisfies
$ u(\cdot,t) \rightarrow u_{\mbox{}_{\!\;\!0}} \!\;\!$ \linebreak
in $ L^{1}_{\tt loc}(\mathbb{R}) $ as $ t \rightarrow 0 $.
It follows from the a priori estimates given in
Section~2 below \linebreak
that all solutions of problem (1.1$a$), (1.1$b$)
are actually globally defined
($ T_{\mbox{}_{\scriptstyle \!\ast}} \!= \infty $), \linebreak
with
$ {\displaystyle
u(\cdot,t) \in
% L^{\infty}_{\tt loc}
% (\;\![\;\!0, \infty \:\! [, L^{p}(\mathbb{R}))
C^{0}(\;\![\;\!0, \infty \:\! [, L^{p}(\mathbb{R}))
} $
for each $ \:\!p \geq p_{\mbox{}_{\!\;\!0}} \!\:\!$ finite.
Given
$ \;\!b \in L^{\infty}(\mathbb{R} \times [\;\!0, \infty\;\![\;\!) $,
what then can be said about the size of supnorm values
$ {\displaystyle
\|\, u(\cdot,t) \,\|_{\mbox{}_{\scriptstyle L^{\infty}(\mathbb{R})}}
} $
for $ t \gg 1 $?

\nl
\mbox{} \vspace{-2.50cm} \\
\par
When
% $ b_{x} \!\equiv \partial b/\partial x \geq 0 $ for all $ x, t $
$ \partial b/\partial x \geq 0 \;\!$ for all $ x \in \mathbb{R}, t \geq 0 $,
it is well known that,
for each $ \:\!p_{\mbox{}_{\!\;\!0}} \!\leq p \leq \;\!\!\infty $,
$ {\displaystyle
\:\!
\|\, u(\cdot,t) \,\|_{\mbox{}_{\scriptstyle L^{p}(\mathbb{R})}}
\!
} $
is monotonically decreasing in $t$,
% decreases monotonically with~$t$,
with \\
\mbox{} \vspace{-0.750cm} \\
\begin{equation}
\tag{1.2}
\|\, u(\cdot,t) \,\|_{\mbox{}_{\scriptstyle L^{\infty}(\mathbb{R})}}
\leq\;\!
K\!\;\!(p_{\mbox{}_{\!\;\!0}}) \,
\|\, u_{\mbox{}_{0}} \;\!
\|_{\mbox{}_{\scriptstyle L^{p_{\mbox{}_{\!\;\!0}}}(\mathbb{R})}} \;\!
t^{\mbox{}^{\scriptstyle \! - \,
\frac{\scriptstyle 1}
     {\scriptstyle \;\!2 \;\! p_{\mbox{}_{\mbox{}_{\!0}}}\!\:\!} }}
\qquad
\forall \;\, t > 0
\qquad
\;\;(\;\! b_x \geq 0 \;\!)\!\!\!\!
\end{equation}
\mbox{} \vspace{-0.250cm} \\
for some constant
$\;\!0 <\!\;\! K\!\;\!(p_{\mbox{}_{\!\;\!0}}) <
\,\! 2^{\mbox{}^{\scriptstyle \!-\;\!1/p_{\mbox{}_{\!\;\!0}} }} \!$
that depends only on $ p_{\mbox{}_{\!\;\!0}} \;\!\!$,
% but not on $b$, $ u_{\mbox{}_{0}} \!\:\!$ or $u$,
see e.g.$\;$\cite{AmickBonaSchonbek1989, %
BrazSchutzZingano2013, EscobedoZuazua1991, %
Porzio2009, Schonbek1986}.
For general $ b(x,t)$, however,
estimating
$ {\displaystyle
\;\!
\|\, u(\cdot,t) \,\|_{\mbox{}_{\scriptstyle L^{\infty}(\mathbb{R})}}
\!
} $
% is not so simple.
is much harder.
% does not seem so simple.
To see why,
let us illustrate with
the important case
$ p_{\mbox{}_{\!\;\!0}} \!= 1 $,
% where one has\footnote{%
% % See e.g.$\,$Theorem 2.1 in Section 2.
% This is reviewed in Theorem 2.1 below.
% } \\
where one has \\
%
% \mbox{} \vspace{-1.075cm} \\
\mbox{} \vspace{-0.750cm} \\
\begin{equation}
\tag{1.3}
\|\, u(\cdot,t) \,\|_{\mbox{}_{\scriptstyle L^{1}(\mathbb{R})}}
\!\;\!\leq\;
\|\, u_{\mbox{}_{0}} \;\!\|_{\mbox{}_{\scriptstyle L^{1}(\mathbb{R})}}
\qquad
\forall \;\, t > 0,
\end{equation}
\mbox{} \vspace{-0.300cm} \\
as recalled in Theorem 2.1 below.
Writing equation (1.1$a$) as \\
\mbox{} \vspace{-0.800cm} \\
\begin{equation}
\tag{1.4}
u_t \;\!+\, b(x,t) \;\!u_x
\;\!=\;
u_{xx} \;\!-\, b_{x}(x,t) \;\!u,
\end{equation}
\mbox{} \vspace{-0.375cm} \\
we observe on the righthand side of (1.4)
that $ |\,u(x,t)\,| $ is pushed to grow
at points $ (x,t) $ where
$ b_{x}(x,t) < 0 $.
If this condition persists long enough,
large values of $ |\,u(x,t)\,| $ might be generated,
particularly at sites where \mbox{$ - \;\!b_{x}(x,t) \gg 1 $}.
% is big.
Now,
% by~(1.3),
because of the constraint (1.3),
any persistent growth in solution size
% will eventually result in the formation
% of long thin structures of the type shown
will eventually create
long thin structures as shown
in Fig.$\,$1,
% which, in turn, are more strongly
% subjected to viscous dissipation,
% which, in turn, are more effectively
% dissipated by the viscous term $ u_{xx} $.
which, in turn, tend to be effectively
% dissipated by the viscous term $ u_{xx} $.
dissipated by viscosity.
The final overall behavior
that ultimately results from such competition
is not immediately clear,
either on physical or mathematical grounds. \\
%
% % If we recall that we are merely assuming
% % $ b \in L^{\infty}(\mathbb{R} \times [\;\!0, \infty\;\![) $
% If we recall that we are not assuming very much
% on the values of $ b_x $
% other than having
% $ b \in L^{\infty}(\mathbb{R} \times [\;\!0, \infty\;\![) $
% (and so, in particular,
% $ - \;\! b_x(x,t) $ may get arbitrarily large)
%

%
% * ---------------------------------------------- *
% *                                                *
% *            Include Figure 1 here               *
% *                                                *
% * ---------------------------------------------- *
%
%
\mbox{} \vspace{-1.500cm} \\
\begin{figure}[hb]
\centering
\includegraphics[width = 11.50cm]{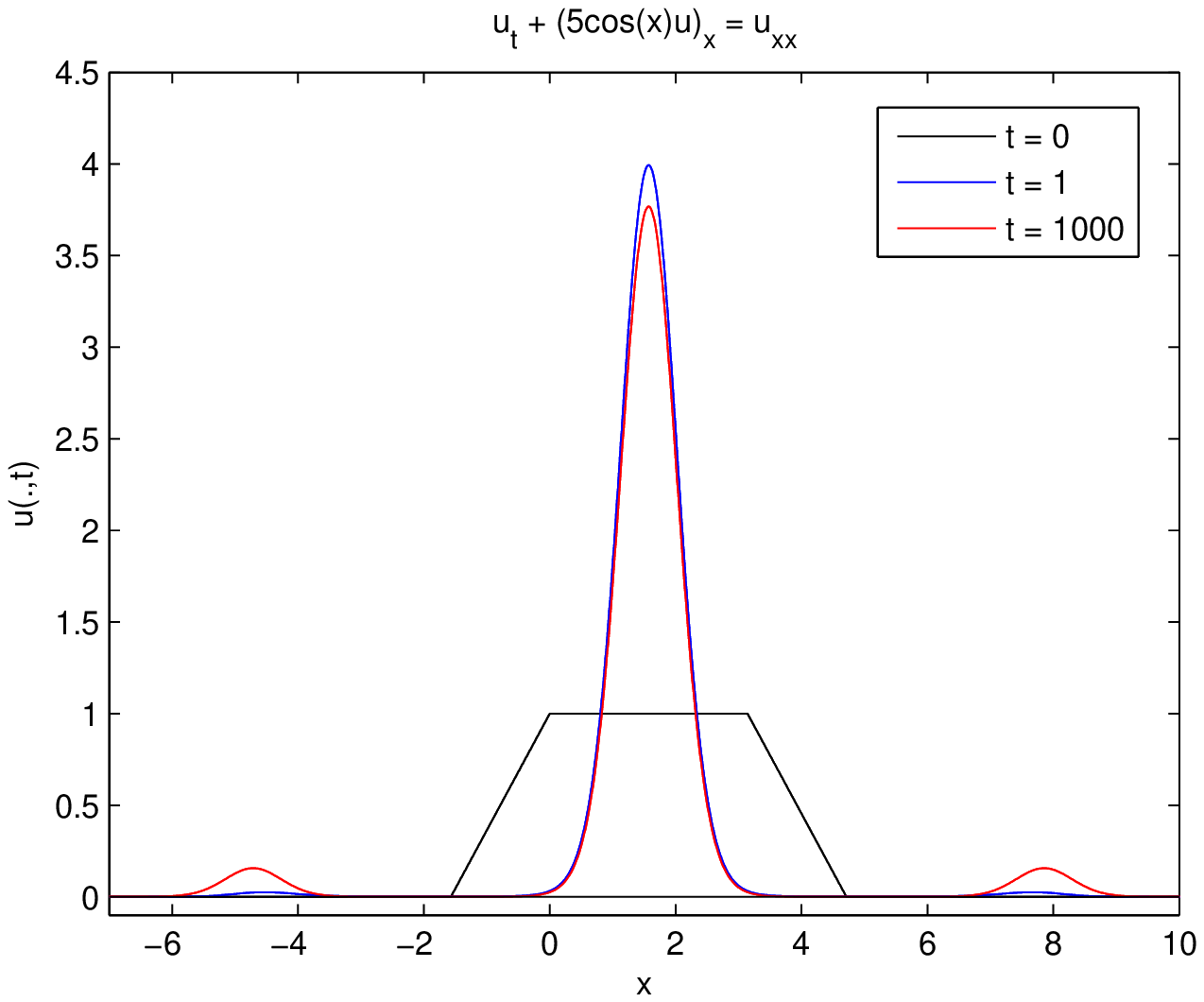}
\end{figure}

\mbox{} \vspace{-1.450cm} \\
\mbox{} \hspace{+1.000cm}
\begin{minipage}[t]{12.50cm}
{\footnotesize
{\bf Fig. 1.}
Solution profiles showing typical growth
in regions with $\;\! b_x \!< 0 $,
% where $ b = 5 \;\!\cos \;\!x $. \linebreak
where \linebreak
%
% \mbox{}$\;$As $ t \rightarrow \infty $,
% solution will slowly spread to the entire line,
% and
% $ \;\!\|\, u(\cdot,t) \,\|_{L^{\infty}(\mathbb{R})} \!\rightarrow 0 $.
%
%
\mbox{} \vspace{-0.690cm} \\
$ b \:\!=\:\! 5\:\!\cos \:\!x $.
After reaching maximum height,
solution starts decaying very slowly \linebreak
\mbox{} \vspace{-0.690cm} \\
due to its spreading
and mass conservation.
(Decay rate is not presently known.)
}
\end{minipage}

\nl
\mbox{} \vspace{-2.500cm} \\
\par
As shown by equation (1.4),
it is not the {magnitude} of $ b(x,t) $
itself
but instead its {\em oscillation\/}
% that plays any role in determining the size of
that is relevant in determining
$ {\displaystyle
\;\!
\|\, u(\cdot,t) \,\|_{\mbox{}_{\scriptstyle L^{\infty}(\mathbb{R})}}
\!\;\!
} $.
Accordingly,
we introduce the quantity $ B(t) $
defined by \\
\mbox{} \vspace{-0.700cm} \\
\begin{equation}
\tag{1.5}
B(t) \,=\,
\frac{\mbox{\small $1$}}{\;\!\mbox{\small $2$}\;\!} \,
\Bigl(\,
\sup_{x \,\in\, \mathbb{R}} \;\!b(x,t)
\;\;\!-\:
\inf_{x \,\in\, \mathbb{R}} \!\;\!b(x,t)
\,\Bigr),
\qquad
t \geq 0,
\end{equation}
\mbox{} \vspace{-0.250cm} \\
which plays a fundamental role in the analysis.
Our main result is now easily stated. \linebreak
\nl
%
% *-----------------------------------------------------------*
% *                                                           *
% *                     Main Theorem                          *
% *                                                           *
% *-----------------------------------------------------------*
%
{\bf Main Theorem.}
{\em
For each $ \;\! p \geq p_{\mbox{}_{\!\;\!0}} $,
we have\/}\footnote{%
In (1.6), (1.11) and other similar expressions
in the text,
it is assumed that
$ 0 \cdot \infty = \infty $.
}
\\
\mbox{} \vspace{-0.950cm} \\
\begin{equation}
\tag{1.6}
\limsup_{t\,\rightarrow\,\infty}\;\!
\|\, u(\cdot,t) \,\|_{\mbox{}_{\scriptstyle L^{\infty}(\mathbb{R})}}
\leq\,
\Bigl(\;\! \frac{\mbox{\small $\;\!3 \;\! \sqrt{\,\!3\;\!} \;$}}
                {\mbox{\small $ 2 \:\! \pi $} } \: p \;\!
\Bigr)^{\mbox{}^{\scriptstyle \!\!
\frac{\scriptstyle 1}{\scriptstyle p} }}
\!\!\cdot\,
{\cal B}^{\mbox{}^{\scriptstyle
\frac{\scriptstyle 1}{\scriptstyle p} }}
\!\!\cdot\:
\limsup_{t\,\rightarrow\,\infty}\;\!
\|\, u(\cdot,t) \,\|_{\mbox{}_{\scriptstyle L^{p}(\mathbb{R})}}
\!,
\end{equation}
\mbox{} \vspace{-0.100cm} \\
{\em
where
}
$ {\displaystyle
{\cal B} =\;\!
\limsup_{t\,\rightarrow \;\!\infty}
\:\!
B(t)
} $.\\
%
% --------------------------------------------------------
%
\nl
In particular,
in the important case
$ p_{\mbox{}_{\!0}} \!= 1 $
considered above,
we obtain, using (1.3), \\
\mbox{} \vspace{-0.675cm} \\
\begin{equation}
\tag{1.7}
\limsup_{t\,\rightarrow\,\infty}\;\!
\|\, u(\cdot,t) \,\|_{\mbox{}_{\scriptstyle L^{\infty}(\mathbb{R})}}
\leq\,
\Bigl(\;\! \frac{\mbox{\small $\;\!3 \;\! \sqrt{\,\!3\;\!} \;$}}
                {\mbox{\small $ 2 \:\! \pi $} } \;\!
\Bigr)
\cdot\,
{\cal B}
\cdot\:
\|\, u_{\mbox{}_{0}} \;\!\|_{\mbox{}_{\scriptstyle L^{1}(\mathbb{R})}}
\!,
\end{equation}
\mbox{} \vspace{-0.210cm} \\
so that
$ u(\cdot,t) $
stays uniformly bounded for all time
in this case.\footnote{%
The constants $(3\sqrt{3}p/(2\pi))^{1/p}$
in (1.6), (1.7) are not optimal;
minimal values are not known.
}
Estimates similar to (1.6)
% should likewise hold for
can be also shown to hold for
% higher dimensional problems as well, e.g.\\
the $n$-dimensional problem \\
\mbox{} \vspace{-0.650cm} \\
\begin{equation}
\tag{1.8}
u_t \,+\;
\mbox{\tt div}\,
(\;\!\mbox{\boldmath $b$}(x,t) \;\!u \;\!)
\:=\:
\Delta \:\!u,
\qquad
u(\cdot,0) \in
L^{p}(\mathbb{R}^{n}) \cap L^{\infty}(\mathbb{R}^{n}),
\end{equation}
\mbox{} \vspace{-0.300cm} \\
but to simplify our discussion
we consider here the case $ n = 1 $ only.
Our derivation of (1.6),
% which improves an unpublished analysis
which improves some unpublished results
by the third author,
%
% given in Section~3
%
% applies to 1-D only but
% is valid in 1-D only but
uses the
% one-dimensional inequality \\
1-D inequality \\
\mbox{} \vspace{-0.900cm} \\
\begin{equation}
\tag{1.9}
\|\; \mbox{v}\;\|_{\mbox{}_{\scriptstyle L^{\infty}(\mathbb{R})}}
\leq\;
C_{\mbox{}_{\!\infty}} \,
\|\; \mbox{v}\;\|_{\mbox{}_{\scriptstyle L^{1}(\mathbb{R})}}
                 ^{\mbox{}^{\scriptstyle 1/3}}
\|\; \mbox{v}_{x}\,\|_{\mbox{}_{\scriptstyle L^{2}(\mathbb{R})}}
                     ^{\mbox{}^{\scriptstyle 2/3}}
\!,
\qquad
\mbox{v} \in L^{1}(\mathbb{R}) \cap H^{1}(\mathbb{R}),
\end{equation}
\mbox{} \vspace{-0.250cm} \\
where
$ C_{\mbox{}_{\!\infty}} \!= (\;\!3/4\;\!)^{\mbox{}^{\scriptstyle 2/3}} \!$,
and can be readily extended to other
problems of interest
like
1-D systems of viscous conservation laws
(\cite{Melo2011}, Ch.$\,$9)
or the more general equation \\
\mbox{} \vspace{-0.750cm} \\
\begin{equation}
\tag{1.10}
u_t \;\!+\, (\;\!b(x,t,u) \;\!u \;\!)_{x}
\;\!=\;
(\;\!a(x,t,u) \;\!u_x \;\!)_{x},
\qquad
a(x,t,u) \geq \mu(t) > 0,
\end{equation}
\mbox{} \vspace{-0.300cm} \\
with bounded values $ b(x,t,u) $,
provided that we assume
$ {\displaystyle
\!\;\!
\int^{\infty} \!\!\!\! \mu(t) \, dt \,=\, \infty
} $:
using a \linebreak
\mbox{} \vspace{-0.575cm} \\
similar argument,
we get the estimate$^{1}$
(\cite{Oliveira2013}, Ch.$\;$2) \\
\mbox{} \vspace{-0.900cm} \\
\begin{equation}
\tag{1.11}
\limsup_{t\,\rightarrow\,\infty}\;\!
\|\, u(\cdot,t) \,\|_{\mbox{}_{\scriptstyle L^{\infty}(\mathbb{R})}}
\leq\,
\Bigl(\;\! \frac{\mbox{\small $\;\!3 \;\! \sqrt{\,\!3\;\!} \;$}}
                {\mbox{\small $ 2 \:\! \pi $} } \: p \;\!
\Bigr)^{\mbox{}^{\scriptstyle \!\!
\frac{\scriptstyle 1}{\scriptstyle p} }}
\!\!\cdot\,
{\cal B}_{\mu}^{\mbox{}^{\scriptstyle \,
\frac{\scriptstyle 1}{\scriptstyle p} }}
\!\!\cdot\:
\limsup_{t\,\rightarrow\,\infty}\;\!
\|\, u(\cdot,t) \,\|_{\mbox{}_{\scriptstyle L^{p}(\mathbb{R})}}
\!,
\end{equation}
\mbox{} \vspace{-0.200cm} \\
for each
$ \;\! p \geq p_{\mbox{}_{\!\;\!0}} \!\;\!$,
where \\
\mbox{} \vspace{-0.825cm} \\
\begin{equation}
\tag{1.12$a$}
{\cal B}_{\mu} \;\!=\;
\limsup_{t\,\rightarrow \;\!\infty}
\;
% B(t)/\mu(t),
\frac{\;\!B(t)\;\!}{\mu(t)},
\end{equation}
\mbox{} \vspace{-0.800cm} \\
\begin{equation}
\tag{1.12$b$}
B(t) \,=\,
\frac{\mbox{\small $1$}}{\;\!\mbox{\small $2$}\;\!} \,
\Bigl(\,
\sup_{x \,\in\, \mathbb{R}} \;\!b(x,t,u(x,t))
\;\;\!-\:
\inf_{x \,\in\, \mathbb{R}} \!\;\!b(x,t,u(x,t))
\,\Bigr).
\end{equation}
\mbox{} \vspace{-0.050cm} \\
More involving applications,
such as problems with
superlinear advection,
% or degenerate diffusion,
% will be described in the future. \\
where solutions may blow up in finite time,
will be described in a sequel to this work. \\
\mbox{} \vspace{-0.150cm} \\
%
%
% ************************************************************
% *                                                          *
% *           Section 2: Basic a priori estimates            *
% *                                                          *
% ************************************************************
%
\par
{\bf \S 2. A priori estimates} \\
\par
This section contains
some preliminary results on
the solutions
% $ \mbox{\boldmath $u$}(\cdot,t) $
of problem~(1.1)
needed later for our derivation of
estimate (1.6),
which is completed in Section~3. \linebreak
($\,\!$Recall that a solution on some given time interval
$ [\;\!0, \mbox{\small $T$}_{\!\ast} [ $,
$ 0 < \mbox{\small $T$}_{\!\ast} \!\leq \infty $,
is a function
$ {\displaystyle
u(\cdot,t) \in
L^{\infty}_{\tt loc}
(\:\![\;\!0, \mbox{\small $T$}_{\!\ast}[, L^{\infty}(\mathbb{R})\:\!)
} $
which is smooth ($C^{2}$ in $x$, $C^{1}$ in $t$)
in $ \mathbb{R} \;\!\times\, ]\;\!0, \mbox{\small $T$}_{\!\ast} \,\![ $
and solves equation (1.1$a$) there,
verifying the initial condition
in the sense of $L^{1}_{\tt loc}(\mathbb{R}) $,
i.e.,
$ {\displaystyle
\|\, u(\cdot,t) - u_{\mbox{}_{\!\,\!0}} \;\!
\|_{\mbox{}_{\scriptstyle L^{1}(\mathbb{K})}}
\!\!\!\;\!\rightarrow 0
\;\!
} $
as $ \,\!t \rightarrow 0 \,\!$
for each compact $ \mathbb{K} \!\;\!\subset \mathbb{R} $.
% For local existence theory, see e.g.$\;$\cite{Serre1999}, Ch.$\;$6.$\,\!$)
% Local existence theory is given e.g.$\;$in \cite{Serre1999}, Ch.$\;$6.$\,\!$)
Local existence theory can be found in e.g.$\;$\cite{Serre1999}, Ch.$\:$6.$\,\!$)
We start with a simple Gronwall-type estimate
for
$\|\, u(\cdot,t) \,\|_{\mbox{}_{\scriptstyle L^{q}(\mathbb{R})}} \!$,
$ p_{\mbox{}_{\!\;\!0}} \!\,\!\leq q < \infty $.
The corresponding result for the supnorm ($q = \infty$)
is more difficult to obtain and will be given
at the end of Section 2,
see Theorem~2.4. \linebreak
\nl
%
% ------------------------------------------------------- %
%                                                         %
%                      Theorem 2.1                        %
%                                                         %
% ------------------------------------------------------- %
%
{\bf Theorem 2.1.}
\textit{%
If
$ {\displaystyle
\,
u(\cdot,t) \in
L^{\infty}_{\tt loc}(\:\![\;\!0, T_{\ast}[, L^{\infty}(\mathbb{R})\:\!)
\:\!
} $
solves problem $\:\!(1.1a)$, $(1.1b)$,
% is a solution to $(1.1a)$, $(1.1b)$,
% solves problem $\:\!(1.1a)$, $(1.1b)$.
then
$ {\displaystyle
u(\cdot,t) \in
C^{0}(\:\![\;\!0, T_{\ast}[, L^{q}(\mathbb{R})\:\!)
\:\!
} $
for each
$ \;\! p_{\mbox{}_{\!\;\!0}} \!\leq q < \infty $,
and
} \\
\mbox{} \vspace{-0.875cm} \\
\begin{equation}
\tag{2.1}
\|\, u(\cdot,t) \,
\|_{\mbox{}_{\scriptstyle L^{q}(\mathbb{R})}}
\leq\;
\|\, u(\cdot,0) \,
\|_{\mbox{}_{\scriptstyle L^{q}(\mathbb{R})}}
\!\cdot\,
\exp \,\Bigl\{\,
% \mbox{\small $ {\displaystyle \frac{\small 1}{\small 2} }$}
\mbox{\small $ {\displaystyle \frac{\small 1}{\small 2} }$}
\, (q-1) \!
\int_{0}^{\mbox{\footnotesize $\:\!t$}} \!\!
% \mbox{\small $ {\displaystyle \frac{\small \;\!B(\tau)^{2}}{\small \mu(\tau)} }$}
B(\tau)^{2}
\; d\tau \,\Bigr\}
\end{equation}
\mbox{} \vspace{-0.250cm} \\
\textit{%
for all $ \,0 < t < T_{\ast} $.
}

%
% -------------------------------------------------------
%
\mbox{} \vspace{-0.150cm} \\
{\small
{\bf Proof.}
The proof is standard,
so we will only sketch the basic steps.
Taking $ S \in C^{1}(\mathbb{R}) $
such that
$ S^{\prime}({\tt v}) \geq 0 $ for all ${\tt v}$,
$ S(0) = 0 $,
$ S({\tt v}) = \mbox{sgn} \;\!({\tt v}) $
for $ | \,{\tt v} \,| \geq 1 $,
let (given $ \delta > 0 $)
$ L_{\delta}({\tt u}) = \int_{0}^{\mbox{\footnotesize ${\tt u}$}}
  \! S({\tt v}/\delta) \, d{\tt v} $,
so that
$ L_{\delta}({\tt u}) \rightarrow |\,{\tt u}\,| $
as $ \delta \rightarrow 0 $, uniformly in $ {\tt u} $.
Let
$ \Phi_{\delta}({\tt u}) = L_{\delta}({\tt u})^{\mbox{\footnotesize $q$}} \!\;\!$.
\linebreak
Given $ \mbox{\footnotesize $R$} > 0 $, $ 0 < \epsilon \leq 1 $,
let $ \zeta_{\mbox{}_{R}}(\cdot) $
be the cut-off function
$ \;\!\zeta_{\mbox{}_{R}}(x) = 0 \;\!$
for $ \;\!|\,x\,| \geq R $,
$ \,\zeta_{\mbox{}_{R}}(x) = $ \linebreak
$ {\displaystyle
\exp\;\!\{\;\!-\:\epsilon \,\sqrt{1 + x^{2}\,}\,\} \;\!-\:
\exp\;\!\{\;\!-\:\epsilon \,\sqrt{1 + \mbox{\footnotesize $R$}^{2}\,}\,\}
\;\!
} $
for $ \;\!|\,x\,| < \mbox{\footnotesize $R$} $.
Multiplying equation (1.1$a$) by \linebreak
$ \Phi_{\delta}^{\prime}(u(x,t)) \cdot \zeta_{\mbox{}_{R}}(x) \;\! $
if $ q \neq 2 $, or $ \,\!u(x,t) \cdot \zeta_{\mbox{}_{R}}(x) \;\! $
if $ q = 2 $,
and
integrating the result on
$ \mathbb{R} \!\;\!\times\!\;\! [\;\!0, t\;\!] $,
we obtain,
letting $ \delta \rightarrow 0 $
and then
$ \mbox{\footnotesize $R$} \rightarrow \infty $,
since
$ {\displaystyle
\;\!
u \in
L^{\infty}(\mathbb{R} \!\;\!\times\!\;\! [\;\!0, t\;\!])
} $: \\
\mbox{} \vspace{-0.750cm} \\
\begin{equation}
\tag{2.2$a$}
{\tt U}_{\epsilon}(t) \,+\, V_{\epsilon}(t)
\;\leq\;
{\tt U}_{\epsilon}(0)
% \|\, u(\cdot,0) \,\|_{\mbox{}_{L^{q}(\mathbb{R})}}^{\,q}
\,+
\int_{0}^{\mbox{\footnotesize $\:\!t$}} \!\!\:\!
G_{\epsilon}(\tau) \, {\tt U}_{\epsilon}(\tau) \: d\tau,
\quad
\;
{\tt U}_{\epsilon}(t) \,=
\int_{\mathbb{R}} \!\;\!
|\, u(x,t) \,|^{\mbox{}^{\mbox{\scriptsize $q$}}} \;\!
w_{\epsilon}(x) \: dx,
% \exp\;\!\{\;\!-\:\epsilon \,\sqrt{1 + x^{2}\,}\,\} \, dx
\end{equation}
\mbox{} \vspace{-0.100cm} \\
where
$ {\displaystyle
\,
w_{\epsilon}(x) \;\!=\;\! \exp\;\!\{\;\!-\:\epsilon \,\sqrt{1 + x^{2}\,}\,\}
} $,
$ {\displaystyle
\;\!
G_{\epsilon}(t) \;\!=\;\!
\frac{1}{2} \, q \;\!(q-1) \;\! B(t)^{2}
+\,\!
\epsilon \;\! 2 \;\! q \,\!\cdot \!\!\!
\sup_{0\,\leq\,\tau\,\leq\,t} \!\!\!
\|\,u(\cdot,t) \,\|_{\mbox{}_{\scriptstyle L^{\infty}(\mathbb{R})}}
} $ \linebreak
\mbox{} \vspace{-0.5750cm} \\
$ {\displaystyle
+\,
\epsilon
} $,
$\;\!$and \\
\mbox{} \vspace{-0.900cm} \\
\begin{equation}
\tag{2.2$b$}
{\tt V}_{\epsilon}(t) \,=\:
\left\{\,
\begin{array}{lll}
\mbox{$ {\displaystyle
\frac{1}{2} \: q \,(q-1) \!
\int_{0}^{\mbox{\footnotesize $\:\!t$}} \!\!\;\!
\int_{\mbox{\scriptsize $\:\!u \neq 0$}} \hspace{-0.500cm}
|\: u(x,\tau) \,|^{\mbox{}^{\scriptstyle \;\!q\;\!-\;\!2}} \;\!
% \Bigl|\, \frac{\partial \;\! u}{\partial \;\!x} \,\Bigr|^{2}
|\:u_{x}(x,\tau)\,|^{\mbox{}^{\scriptstyle 2}} \;\!
w_{\epsilon}(x) \; dx \, d\tau
} $},
& \mbox{} & \mbox{if }\; q \neq 2, \\
\mbox{} \vspace{-0.250cm} \\
\mbox{$ {\displaystyle
% \frac{1}{2} \: q \,(q-1)
\int_{0}^{\mbox{\footnotesize $\:\!t$}} \!\!\;\!
\int_{\mathbb{R}} \:\!
% \Bigl|\, \frac{\partial \;\! u}{\partial \;\!x} \,\Bigr|^{2}
|\:u_{x}(x,\tau)\,|^{\mbox{}^{\scriptstyle 2}} \;\!
w_{\epsilon}(x) \; dx \, d\tau
} $},
& \mbox{} & \mbox{if }\; q = 2.
\end{array}
\right.
\end{equation}
\mbox{} \vspace{-0.050cm} \\
By Gronwall's lemma,
(2.2) gives
$ \,
{\tt U}_{\epsilon}(t) \leq\;\!
{\tt U}_{\epsilon}(0) \cdot
\exp\,\big\{ \!\;\!\int_{0}^{\mbox{\footnotesize $\:\!t$}}
\!\;\!G_{\epsilon}(\tau) \, d\tau \,\!\bigr\}
$,
from which we obtain \linebreak
(2.1) by simply letting $ \epsilon \rightarrow 0 $.
This shows, in particular,
that
$ {\displaystyle
\;\!
u(\cdot,t) \in
L^{\infty}_{\tt loc}(\:\![\;\!0, T_{\ast}[, L^{q}(\mathbb{R})\:\!)
} $ \linebreak
if $ p_{\mbox{}_{\!\;\!0}} \!\leq q < \infty $.
Now,
to get
$ {\displaystyle
u(\cdot,t) \in
C^{0}(\:\![\;\!0, T_{\ast}[, L^{q}(\mathbb{R})\:\!)
} $,
it is sufficient to show that,
given $\;\!\varepsilon > 0 \;\!$
and
$ \;\! 0 < \mbox{\footnotesize $T$} < \mbox{\footnotesize $T$}_{\ast} $
arbitrary,
we can find
% $ \;\!\mbox{\footnotesize $R$} \gg 1 \;\!$
$ {\displaystyle
\;\!\mbox{\footnotesize $R$} \;\!=\;\!
\mbox{\footnotesize $R$}(\varepsilon,\mbox{\footnotesize $T$})
\gg 1 \;\!
} $
large enough so that
we have
$ {\displaystyle
\;\!
\|\, u(\cdot,t) \,
\|_{\mbox{}_{\scriptstyle L^{q}(\;\!|\,x\,| \,>\, R\;\!)}}
\!\!\;\!< \varepsilon
\:\!
} $
for any
$ \;\!0 \leq t \leq \mbox{\footnotesize $T$} $.
Taking
$ \psi \in C^{2}(\mathbb{R}) $
with $ 0 \leq \psi \leq 1 $
and
$ \psi(x) = 0 $ for all $ x \leq 0 $,
$ \psi(x) = 1 $ for all $ x \geq 1 $,
let
% (for $ \mbox{\footnotesize $R$} > 1 $,
% $ \mbox{\footnotesize $M$} > 1 $)
$ \Psi_{\scriptstyle \!\:\!R,\,M} \in C^{2}(\mathbb{R}) $
be the cut-off function
given by
$ \;\!\Psi_{\scriptstyle \!\:\!R,\,M} (x) = 0 \;\!$
  if $ |\,x\,| \leq \mbox{\footnotesize $R$} - 1 $,
$ \;\!\Psi_{\scriptstyle \!\:\!R,\,M} (x) =
  \psi(\;\!|\,x\,| - \mbox{\footnotesize $R$} + 1 ) \;\!$
  if $ \mbox{\footnotesize $R$} - 1 < |\,x\,| < \mbox{\footnotesize $R$} $,
and \linebreak
$ \;\!\Psi_{\scriptstyle \!\:\!R,\,M} (x) = 1 \;\!$
  if $ \mbox{\footnotesize $R$} \leq |\,x\,| \leq
       \mbox{\footnotesize $R$} + \mbox{\footnotesize $M$} $,
$ \;\!\Psi_{\scriptstyle \!\:\!R,\,M} (x) =
  \psi(\;\!\mbox{\footnotesize $R$} + \mbox{\footnotesize $M$} + 1 - |\,x\,| \;\!) \;\!$
  if $ \mbox{\footnotesize $R$} + \mbox{\footnotesize $M$} < |\,x\,|
       < \mbox{\footnotesize $R$} + \mbox{\footnotesize $M$} + 1 $,
$ \;\!\Psi_{\scriptstyle \!\:\!R,\,M} (x) = 0 \;\!$
  if $ \;\!|\,x\,| \geq \mbox{\footnotesize $R$} + \mbox{\footnotesize $M$} + 1 $,
$\;\!$where
$ \mbox{\footnotesize $R$} > 1 $, $ \mbox{\footnotesize $M$} > 0 $
are given.
Multiplying
(1.1$a$) by
$ \Phi_{\delta}^{\prime}(u(x,t)) \cdot
  \Psi_{\scriptstyle \!\:\!R,\,M} (x)  \;\! $
if $ q \neq 2 $,
or $ \,\!u(x,t) \cdot \Psi_{\scriptstyle \!\:\!R,\,M} (x)  \;\! $
if $ q = 2 $,
and
integrating the result on
$ \mathbb{R} \!\;\!\times\!\;\! [\;\!0, t\;\!] $,
$ 0 < t \leq \mbox{\footnotesize $T$} $,
we obtain,
as in (2.2),
by letting
$ \delta \rightarrow 0 $,
$ \mbox{\footnotesize $M$} \rightarrow \infty $,
that
$ {\displaystyle
\|\, u(\cdot,t) \,
\|_{\mbox{}_{\scriptstyle L^{q}(\;\!|\,x\,| \,>\, R\;\!)}}
\!\!\;\!< \varepsilon/2 \,+\,
\|\, u(\cdot,0) \,
\|_{\mbox{}_{\scriptstyle L^{q}(\;\!|\,x\,| \,>\, R-1\;\!)}}
\:\!
} $
for all
$ 0 \leq t \leq \mbox{\footnotesize $T$} $,
provided that we take $ R > 1 $ sufficiently large.
This gives the continuity result,
and the proof is complete.
\hfill $ \Box$ \\
}
%
% --------------------------------------- End of proof
%                                        for Theorem 2.1

%
\mbox{} \vspace{-1.350cm} \\
\par
An important by-product of the proof above is that
we have
(letting $ \epsilon \rightarrow 0 $ in (2.2),
and using (2.1)),
for each
$ {\displaystyle
\;\!
0 < \mbox{\footnotesize $T$} < \mbox{\footnotesize $T$}_{\ast}
} $
and
$ \;\!q \geq \max\;\!\{\;\! p_{\mbox{}_{\!\;\!0}}, 2 \;\!\} $, \\
\mbox{} \vspace{-0.700cm} \\
\begin{equation}
\tag{2.3}
\int_{0}^{\mbox{\scriptsize $\;\!T$}} \!\!\!
\int_{\mathbb{R}} \,
|\, u(x,\tau) \,|^{\mbox{}^{\scriptstyle q-2}} \,
% \Bigl|\, \frac{\partial \;\! u}{\partial \;\!x} \,\Bigr|^{2}
|\, u_x(x,\tau) \,|^{\mbox{}^{\scriptstyle 2}}
\;\! dx \, d\tau
\,< \infty.
\end{equation}
\mbox{} \vspace{-0.150cm} \\
Therefore,
if we repeat the steps above leading to (2.2),
we obtain
(letting $ \delta \rightarrow 0 $,
$ \mbox{\small $R$} \rightarrow \infty $,
$ \epsilon \rightarrow 0 $,
in this order,
taking (2.1), (2.3) into account)
the identity \\
\mbox{} \vspace{-0.150cm} \\
\mbox{} \hspace{1.00cm}
$ {\displaystyle
\|\, u(\cdot,t) \,
\|_{\mbox{}_{\scriptstyle L^{q}(\mathbb{R})}}
  ^{\mbox{}^{\scriptstyle \;\!q}}
+\;
q \, (q-1) \!
\int_{0}^{\mbox{\footnotesize $\:\!t$}} \!\!\;\!
\int_{\mathbb{R}}
|\, u(x,\tau) \,|^{\mbox{}^{\scriptstyle q - 2}} \,
% \Bigl|\, \frac{\partial \:\! u}{\partial \:\! x} \,\Bigr|^{2}
|\, u_x(x,\tau) \,|^{\mbox{}^{\scriptstyle 2}}
\;\! dx \, d\tau
\;=
} $ \\
\mbox{} \vspace{-0.650cm} \\
\mbox{} \hfill (2.4) \\
\mbox{} \vspace{-0.425cm} \\
\mbox{} \hspace{-0.350cm}
$ {\displaystyle
\mbox{}\;\,=\;\,
\|\, u(\cdot,0) \,
\|_{\mbox{}_{\scriptstyle L^{q}(\mathbb{R})}}
  ^{\mbox{}^{\scriptstyle \;\!q}}
+\;
q \, (q-1) \!
\int_{0}^{\mbox{\footnotesize $\:\!t$}} \!\!
\int_{\mathbb{R}}
\bigl(\;\! b(x,\tau) - \beta(\tau) \:\!\bigr) \,
|\, u(x,\tau) \,|^{\mbox{}^{\scriptstyle q - 2}}
\hspace{-0.300cm}
u(x,\tau) \,
% \frac{\partial \:\! u}{\partial \:\! x}
u_x(x,\tau)
\; dx \, d\tau
} $ \\
\mbox{} \vspace{+0.100cm} \\
for every $ \;\!0 < t < \mbox{\small $T$}_{\!\ast} $
and
$ \;\!\max\;\!\{\;\! p_{\mbox{}_{\!\;\!0}}, 2 \;\!\} \leq q < \infty $,
where \\
\mbox{} \vspace{-0.700cm} \\
\begin{equation}
\tag{2.5}
\beta(t)
\,=\,
\frac{\mbox{\small $1$}}{\;\!\mbox{\small $2$}\;\!} \,
\Bigl(\,
\sup_{x \,\in\, \mathbb{R}} \;\!b(x,t)
\;\;\!+\:
\inf_{x \,\in\, \mathbb{R}} \!\;\!b(x,t)
\,\Bigr),
\qquad
t \geq 0.
\end{equation}
\mbox{} \vspace{-0.250cm} \\
The core of the difficulty in the analysis of (1.1)
is apparent here:
under the sole assumption that
$ b $ is bounded,
it is not much clear
how one should
% deal with
go about
the last term
% on the righthand side of (2.4),
% other than using
% elementary Gronwall-type arguments.
in (2.4) in order to get more than (2.1) above.
Actually,
it will be convenient to consider (2.4)
in the (equivalent) differential form,
i.e., \\
\mbox{} \vspace{-0.110cm} \\
\mbox{} \hspace{1.50cm}
$ {\displaystyle
\frac{d}{d \:\!t} \:
\|\, u(\cdot,t) \,
\|_{\mbox{}_{\scriptstyle L^{q}(\mathbb{R})}}
  ^{\mbox{}^{\scriptstyle \;\!q}}
+\;
q \, (q-1) \!
\int_{\mathbb{R}}
|\, u(x,t) \,|^{\mbox{}^{\scriptstyle q - 2}} \,
% \Bigl|\, \frac{\partial \:\! u}{\partial \:\! x} \,\Bigr|^{2}
|\, u_x(x,t) \,|^{\mbox{}^{\scriptstyle 2}}
\;\! dx
\;=
} $ \\
\mbox{} \vspace{-0.550cm} \\
\mbox{} \hfill (2.6) \\
\mbox{} \vspace{-0.550cm} \\
%
% \mbox{} \hspace{+1.200cm}
\mbox{} \hspace{+2.500cm}
$ {\displaystyle
\mbox{}\;\,=\;
q \, (q-1) \!
\int_{\mathbb{R}}
\bigl(\;\! b(x,t) - \beta(t) \:\!\bigr) \,
|\, u(x,t) \,|^{\mbox{}^{\scriptstyle q - 2}}
\hspace{-0.300cm}
u(x,t) \,
% \frac{\partial \:\! u}{\partial \:\! x}
u_x(x,t)
\; dx
% \quad
% \;\;
% \mbox{a.e.}\;\,
% t \in \:\![\;\!0, \mbox{\small $T$}_{\!\ast} [.
} $ \\
\mbox{} \vspace{+0.100cm} \\
%
% for a.e.
% $ t \in \:\![\;\!0, \mbox{\small $T$}_{\!\ast} [ $.
%
for all
$ {\displaystyle
t \in \:\![\;\!0, \mbox{\small $T$}_{\!\ast} [ \,\setminus\;\! E_{q}
} $,
where $ E_{q} \!\;\!\subset [\;\!0, \mbox{\small $T$}_{\!\ast}[ $
has zero measure.
We then readily obtain,
using (1.9) and the one-dimensional
Nash inequality \cite{Nash1958} \\
\mbox{} \vspace{-0.800cm} \\
\begin{equation}
\tag{2.7}
\|\: {\tt v} \:\|_{\mbox{}_{\scriptstyle L^{2}(\mathbb{R})}}
\leq\,
C_{\mbox{}_{2}} \,
\|\: {\tt v} \:\|_{\mbox{}_{\scriptstyle L^{1}(\mathbb{R})}}
                 ^{\mbox{}^{\scriptstyle \;\!2/3}}
\;\!
\|\: {\tt v}_{x} \,\|_{\mbox{}_{\scriptstyle L^{2}(\mathbb{R})}}
                     ^{\mbox{}^{\scriptstyle \:\!1/3}}
\!,
\qquad
C_{\mbox{}_{2}}
=\,
\Bigl(\;\!
\mbox{\small $ {\displaystyle
\frac{\small \;\!3 \;\!\sqrt{\:\!3\;}\,}{\small 4 \:\! \pi} }$}
\;\!\Bigr)^{\mbox{}^{\scriptstyle \!\!\!\,\! 1/3}}
\!\!\!,
\end{equation}
\mbox{} \vspace{-0.150cm} \\
where the value given above for $C_{\mbox{}_{2}}$
is optimal \cite{CarlenLoss1993},
the following result: \\
\mbox{} \vspace{-0.050cm} \\
%
% ------------------------------------------------------- %
%                                                         %
%                      Theorem 2.2                        %
%                                                         %
% ------------------------------------------------------- %
%
{\bf Theorem 2.2.}
\textit{%
$\!\:\!$Let $ \:\!q \geq 2 \:\!p_{\mbox{}_{\!\;\!0}} $.
$\!$If
$ {\displaystyle
\;\!
\hat{t} \in
\:\![\;\!0, \mbox{\small $T$}_{\!\ast} [ \,\setminus\;\! E_{q}
\!\;\!
} $
is such that
$ {\displaystyle
\;\!
\mbox{\footnotesize $ {\displaystyle \frac{d}{d\:\!t} }$} \,
\|\, u(\cdot,t) \,
\|_{\mbox{}_{\scriptstyle L^{q}(\mathbb{R})}}
  ^{\mbox{}^{\scriptstyle \;\!q}}
{\mbox{}_{\bigr|}}_{\mbox{}_{\mbox{\footnotesize $t = \:\!\hat{t}$}}}
\hspace{-0.700cm}
\geq\: 0
} $, \linebreak
\mbox{} \vspace{-0.550cm} \\
then
}
\mbox{} \vspace{-0.550cm} \\
\begin{equation}
\tag{2.8$a$}
\|\, u(\cdot,\hat{t}\:\!) \,
\|_{\mbox{}_{\scriptstyle L^{q}(\mathbb{R})}}
\;\!\leq\,
\Bigl(\;\!
\mbox{\small $ {\displaystyle \frac{\;\!q\;\!}{2} }$} \,
C_{\mbox{}_{\!\;\!2}}^{\mbox{}^{\scriptstyle \;\!3}}
\:\!\Bigr)^{\mbox{}^{\scriptstyle \!\!\!\;\!1/q}}
\!
B(\:\!\hat{t}\:\!)^{\mbox{}^{\scriptstyle \!\:\!1/q}}
\,
\|\, u(\cdot,\hat{t}\:\!) \,
\|_{\mbox{}_{\scriptstyle L^{q/2}(\mathbb{R})}}
\end{equation}
\mbox{} \vspace{-0.450cm} \\
{\em and} \\
\mbox{} \vspace{-1.150cm} \\
\begin{equation}
\tag{2.8$b$}
\|\, u(\cdot,\hat{t}\:\!) \,
\|_{\mbox{}_{\scriptstyle L^{\infty}(\mathbb{R})}}
\;\!\leq\,
\Bigl(\;\!
\mbox{\small $ {\displaystyle \frac{\;\!q\;\!}{2} }$} \,
C_{\mbox{}_{\!\;\!2}} \;\! C_{\mbox{}_{\!\infty}}
\Bigr)^{\mbox{}^{\scriptstyle \!\!2/q}}
B(\:\!\hat{t}\:\!)^{\mbox{}^{\scriptstyle \!\;\!2/q}}
\,
\|\, u(\cdot,\hat{t}\:\!) \,
\|_{\mbox{}_{\scriptstyle L^{q/2}(\mathbb{R})}}
\!\,\!.
\end{equation}
\mbox{} \vspace{+0.150cm} \\
%
%
% -------------------------------------------------------
%
\mbox{} \vspace{-0.600cm} \\
{\small
{\bf Proof.}
Consider (2.8$a$) first.
From (1.5), (2.5) and (2.6),
we have \\
\mbox{} \vspace{-0.675cm} \\
\begin{equation}
\notag
\int_{\mathbb{R}} \!\;\!
|\, u(x,\hat{t}\:\!) \,|^{\mbox{}^{\scriptstyle q - 2}}
\;\!
% \Bigl|\;\! \frac{\partial \:\! u}{\partial \:\!x} \;\!\Bigr|^{2}
|\, u_x(x,\hat{t}\:\!) \,|^{\mbox{}^{\scriptstyle 2}}
\;\!dx
\;\leq\;
B(\:\!\hat{t}\:\!) \!
\int_{\mathbb{R}} \!\;\!
|\, u(x,\hat{t}\:\!) \,|^{\mbox{}^{\scriptstyle q - 1}}
\;\!
|\;\! u_x(x,\hat{t}\:\!) \,|
\;dx.
\end{equation}
\mbox{} \vspace{-0.100cm} \\
This gives \\
\mbox{} \vspace{-0.800cm} \\
\begin{equation}
\notag
\int_{\mathbb{R}} \!\;\!
|\, u(x,\hat{t}\:\!) \,|^{\mbox{}^{\scriptstyle q - 2}}
\,
% \Bigl|\;\! \frac{\partial \:\! u}{\partial \:\!x} \;\!\Bigr|^{2}
|\, u_x(x,\hat{t}\:\!) \,|^{\mbox{}^{\scriptstyle 2}}
\;\!dx
\;\leq\;
B(\:\!\hat{t}\:\!)^{\mbox{}^{\scriptstyle 2}}
\;\!
\|\, u(\cdot,\hat{t}\:\!) \,
\|_{\mbox{}_{\scriptstyle L^{q}(\mathbb{R})}}
  ^{\mbox{}^{\scriptstyle \;\!q}}
\!,
\end{equation}
\mbox{} \vspace{-0.200cm} \\
or,
in terms of
$ {\displaystyle
\;\!
\hat{v} \in L^{1}(\mathbb{R}) \cap L^{\infty}(\mathbb{R})
} $
defined by
$ {\displaystyle
\;\!
\hat{v}(x) =
|\,u(x,\hat{t}\:\!)\,|^{\mbox{}^{\scriptstyle \;\!q/2}}
\!\;\!
} $
if $ q > 2 $,
$ {\displaystyle
\;\!
\hat{v}(x) = u(x,\hat{t}\:\!)
\;\!
} $
if $ q = 2 $, \\
\mbox{} \vspace{-0.950cm} \\
\begin{equation}
\notag
\|\: \hat{v}_{x} \,
\|_{\mbox{}_{\scriptstyle L^{2}(\mathbb{R})}}
\,\leq\;
\frac{\;\!q\;\!}{2} \,
B(\:\!\hat{t}\:\!) \;
\|\: \hat{v} \:
\|_{\mbox{}_{\scriptstyle L^{2}(\mathbb{R})}}\!\:\!.
\end{equation}
\mbox{} \vspace{-0.100cm} \\
Using (2.7),
we then get
$ {\displaystyle
\;\!
\|\: \hat{v} \:
\|_{\mbox{}_{\scriptstyle L^{2}(\mathbb{R})}}
  ^{\mbox{}^{\scriptstyle \;\!2}}
\leq\,
\frac{\;\!q\;\!}{2} \:
C_{\mbox{}_{\!2}}^{\mbox{}^{\scriptstyle \;\!3}}
\;\!
B(\:\!\hat{t}\:\!)  \:
\|\: \hat{v} \:
\|_{\mbox{}_{\scriptstyle L^{1}(\mathbb{R})}}
  ^{\mbox{}^{\scriptstyle \;\!2}}
\!\:\!
} $,
which is equivalent to (2.8$a$). \linebreak
\mbox{} \vspace{-0.525cm} \\
Similarly,
(2.8$b$) can be obtained,
using (1.9).
\mbox{}
\hfill $ \Box$ \\
}
%
% --------------------------------------- End of proof
%                                        for Theorem 2.2
%
\mbox{} \vspace{-0.950cm} \\
% xxxxxxxxxxxxxxxxxxxxxxxxxxxx
\par
Thus, we can use (2.8) when
$ {\displaystyle
\;\!
\|\,u(\cdot,t) \,
\|_{\mbox{}_{\scriptstyle L^{q}(\mathbb{R})}}
\!
} $
is not decreasing.
If it is decreasing,
%
% On the other hand,
% at $t$ values where
% $ {\displaystyle
% \;\!
% \|\, u(\cdot,t) \,
% \|_{\mbox{}_{\scriptstyle L^{q}(\mathbb{R})}}
% \!
% } $
% happens to be decreasing, \linebreak
% %
% \mbox{} \vspace{-0.570cm} \\
% %
(2.6) becomes useless but
at least
% we know something already
% about
we know in such case that
$ {\displaystyle
\;\!
\|\, u(\cdot,t) \,
\|_{\mbox{}_{\scriptstyle L^{q}(\mathbb{R})}}
} $
% in this case that can also be used.
% % The next result shows how this can be done.
% The next result gives an effective way to proceed.
is not increasing,
which should be useful too.
Different values of $q$
have
different scenarios,
% which have to be pieced together in some way.
which we will have to piece together in some way.
The next result shows us just how.
To this end,
it is convenient
to introduce the quantities
$ \mathbb{B}(t_0\:\!; t) $,
$ \mathbb{U}_{\!\;\!p}(t_0\:\!; t) $
defined by \\
\mbox{} \vspace{-0.570cm} \\
\begin{equation}
\tag{2.9}
\mathbb{B}(t_0\:\!; t)
\;=\;\;\!
\sup\:\Bigl\{\;\!
B(\tau)\!\;\!:
\; t_0 \!\leq \tau \leq t \;\Bigr\},
\end{equation}
\mbox{} \vspace{-0.825cm} \\
\begin{equation}
\tag{2.10}
\mbox{} \;\;
\mathbb{U}_{p}(t_0\:\!; t)
\;=\;\;\!
\sup\:\Bigl\{\;\!
\|\, u(\cdot,\tau) \,
\|_{\mbox{}_{\scriptstyle L^{p}(\mathbb{R})}}
\!\!\;\!:
\; t_0 \!\leq \tau \leq t \;\Bigr\},
\end{equation}
\mbox{} \vspace{-0.175cm} \\
given
$ \;\!p \geq p_{\mbox{}_{\!\;\!0}} \!\;\!$,
$ \;\! 0 \leq t_0 \!\;\!\leq t < \mbox{\small $T$}_{\!\!\;\!\ast} $
arbitrary. \\
\nl
%
% ------------------------------------------------------- %
%                                                         %
%                      Theorem 2.3                        %
%                                                         %
% ------------------------------------------------------- %
%
{\bf Theorem 2.3.}
\textit{%
Let $\;\! q \geq 2 \:\!p_{\mbox{}_{\!\;\!0}} $.
For each
$ \,0 \leq t_0 \!\;\!< \mbox{\small $T$}_{\!\!\;\!\ast} $,
we have
} \\
\mbox{} \vspace{-0.900cm} \\
\begin{equation}
\tag{2.11}
\mathbb{U}_{\!\;\!q}(t_0\:\!; t)
\;\leq\;\;\!
\max\,\biggl\{\;\!
\|\, u(\cdot,t_0) \,
\|_{\mbox{}_{\scriptstyle L^{q}(\mathbb{R})}};
\,
\Bigl(\;\! \frac{\;\!q\;\!}{2} \,
C_{\mbox{}_{\!2}}^{\mbox{}^{\scriptstyle \;\!3}}
\,\!\Bigr)^{\mbox{}^{\scriptstyle \!\!\!
\frac{\scriptstyle 1}{\scriptstyle q} }}
\mathbb{B}(t_0\:\!; t)^{\mbox{}^{\scriptstyle \!
\frac{\scriptstyle 1}{\scriptstyle q} }}
\;\!
\mathbb{U}_{\mbox{}_{\scriptstyle \!
\frac{\scriptstyle q}{\scriptscriptstyle 2}}}\!(t_0\:\!; t)
\,\biggr\}
\end{equation}
\mbox{} \vspace{-0.200cm} \\
\textit{%
for all}
$\;\! t_0 \!\;\! \leq t < \mbox{\small $T$}_{\!\!\;\!\ast} $. \\
%
%
% -------------------------------------------------------
%
\mbox{} \vspace{-0.500cm} \\
{\small
{\bf Proof.}
Set
$ {\displaystyle
\;\!
\lambda_{q}(t)
\;\!=\;\!
\Bigl(\;\! \frac{\;\!q\;\!}{2} \,
C_{\mbox{}_{\!2}}^{\mbox{}^{\scriptstyle \;\!3}}
\,\!\Bigr)^{\mbox{}^{\scriptstyle \!\!\!
\frac{\scriptstyle 1}{\scriptstyle q} }}
\mathbb{B}(t_0\:\!; t)^{\mbox{}^{\scriptstyle \!
\frac{\scriptstyle 1}{\scriptstyle q} }}
\;\!
\mathbb{U}_{\mbox{}_{\scriptstyle \!
\frac{\scriptstyle q}{\scriptscriptstyle 2}}}\!(t_0\:\!; t)
} $.
% The proof includes three cases: \\
% We have three cases to consider: \\
There are three cases to consider: \\
\mbox{} \vspace{-0.300cm} \\
%
%
% ---------------------------------------- Case I
%
%
{\tt Case I:}
$ {\displaystyle
\|\, u(\cdot,\tau) \,
\|_{\mbox{}_{\scriptstyle L^{q}(\mathbb{R})}}
\!> \lambda_{q}(t)
\;\!
} $
for all $\;\! t_0 \!\;\!\leq \tau \leq t $.
$\;\!$By (2.8$a$), Theorem 2.2,
we must then have
$ {\displaystyle
\;\!
d/d\tau \;\!
\|\, u(\cdot,\tau) \,
\|_{\mbox{}_{\scriptstyle L^{q}(\mathbb{R})}}
  ^{{\scriptstyle \;\!q}}
\!\!\;\!< 0
\;\!
} $
for all
$ \;\!\tau \!\;\!\in \!\;\![\;\!t_0, t \;\!] \:\!\setminus\,\! E_{q} $,
so that
$ {\displaystyle
\;\!
\|\, u(\cdot,\tau) \,
\|_{\mbox{}_{\scriptstyle L^{q}(\mathbb{R})}}
\!
} $
is monotonically decreasing in
$ [\;\!t_0, t \;\!] $.
In particular,
$ {\displaystyle
\;\!
\mathbb{U}_{\!\;\!q}(t_0\:\!; t)
\!\;\!=\;\!
\|\, u(\cdot,t_0) \,
\|_{\mbox{}_{\scriptstyle L^{q}(\mathbb{R})}}
\!
} $
in this case,
and (2.11) holds. \\
%
% ---------------------------------------- Case II
%
\mbox{} \vspace{-0.300cm} \\
{\tt Case II:}
$ {\displaystyle
\|\, u(\cdot,t_0) \,
\|_{\mbox{}_{\scriptstyle L^{q}(\mathbb{R})}}
\!> \lambda_{q}(t)
\;\!
} $
and
$ {\displaystyle
\;\!
\|\, u(\cdot,t_1) \,
\|_{\mbox{}_{\scriptstyle L^{q}(\mathbb{R})}}
\!\leq \lambda_{q}(t)
\;\!
} $
for some
$ t_1 \!\in\; ]\;\!t_0, t\;\!] $. \\
\mbox{} \vspace{-0.500cm} \\
In this case,
let
$ t_2 \!\in\; ]\;\!t_0, t\;\!] $
be such that
we have
$ {\displaystyle
\;\!
\|\, u(\cdot,\tau) \,
\|_{\mbox{}_{\scriptstyle L^{q}(\mathbb{R})}}
\!>\!\;\! \lambda_{q}(t)
\;\!
} $
for all
$ t_0 \!\leq \tau < t_2 $,
while
$ {\displaystyle
\;\!
\|\, u(\cdot,t_2) \,
\|_{\mbox{}_{\scriptstyle L^{q}(\mathbb{R})}}
\!= \lambda_{q}(t)
} $.
We claim that
$ {\displaystyle
\;\!
\|\, u(\cdot,\tau) \,
\|_{\mbox{}_{\scriptstyle L^{q}(\mathbb{R})}}
\!\leq\!\;\! \lambda_{q}(t)
\;\!
} $
for every
$ t_2 \!\leq \tau \leq t $: \linebreak
in fact,
if this were not true,
we could then find
$ t_3, t_4 $
with
$ t_2 \!\leq t_3 \!< t_4 \!\leq t $
such that
$ {\displaystyle
\|\, u(\cdot,\tau) \,
\|_{\mbox{}_{\scriptstyle L^{q}(\mathbb{R})}}
\!>\!\;\! \lambda_{q}(t)
\;\!
} $
for all
$ t_3 \!< \tau \leq t_4 $,
$ {\displaystyle
\;\!
\|\, u(\cdot,t_3) \,
\|_{\mbox{}_{\scriptstyle L^{q}(\mathbb{R})}}
\!= \lambda_{q}(t)
} $.
By (2.8$a$), Theorem 2.2, \linebreak
this would require
$ {\displaystyle
\,
d/d\tau \;\!
\|\, u(\cdot,\tau) \,
\|_{\mbox{}_{\scriptstyle L^{q}(\mathbb{R})}}
  ^{{\scriptstyle \;\!q}}
\!\! < 0
\;\!
} $
for all
$ \;\!\tau \!\;\!\in \;]\;\!t_3, t_4 \,\!] \:\!\setminus\,\! E_{q} $,
so that
$ {\displaystyle
\;\!
\|\, u(\cdot,\tau) \,
\|_{\mbox{}_{\scriptstyle L^{q}(\mathbb{R})}}
\!
} $  \linebreak
could not increase anywhere on
$ \;\![\;\!t_3, t_4 \,\!] $.
This contradicts
$ {\displaystyle
\;\!
\|\, u(\cdot,t_3) \,
\|_{\mbox{}_{\scriptstyle L^{q}(\mathbb{R})}}
\!<\;\!
\|\, u(\cdot,t_4) \,
\|_{\mbox{}_{\scriptstyle L^{q}(\mathbb{R})}}
\!\;\!
} $,
and so
we have
$ {\displaystyle
\;\!
\|\, u(\cdot,\tau) \,
\|_{\mbox{}_{\scriptstyle L^{q}(\mathbb{R})}}
\!\leq\!\;\! \lambda_{q}(t)
\;\!
} $
for every
$ t_2 \!\leq \tau \leq t $,
as claimed.
On the other hand,
by (2.8$a$),
$ {\displaystyle
\|\, u(\cdot,\tau) \,
\|_{\mbox{}_{\scriptstyle L^{q}(\mathbb{R})}}
\!
} $
has to be monotonically decreasing on
$ \;\![\;\!t_0, t_2 \,\!] $,
just as in {\tt Case}~{\tt I}.
Therefore,
we have
$ {\displaystyle
\;\!
\mathbb{U}_{\!\;\!q}(t_0\:\!; t)
\!\;\!=\;\!
\|\, u(\cdot,t_0) \,
\|_{\mbox{}_{\scriptstyle L^{q}(\mathbb{R})}}
\!
} $
in this case again,
which shows (2.11). \\
%
%
% ---------------------------------------- Case III
%
\mbox{} \vspace{-0.300cm} \\
{\tt Case III:}
$ {\displaystyle
\|\, u(\cdot,t_0) \,
\|_{\mbox{}_{\scriptstyle L^{q}(\mathbb{R})}}
\!\leq \lambda_{q}(t)
} $.
This gives
$ {\displaystyle
\;\!
\|\, u(\cdot,\tau) \,
\|_{\mbox{}_{\scriptstyle L^{q}(\mathbb{R})}}
\!\leq\!\;\! \lambda_{q}(t)
\;\!
} $
for every
$ t_0 \!\leq \tau \leq t $,
by repeating the argument used
on the interval
$\;\![\;\!t_2, t \;\!] \;\!$
in {\tt Case II} above.
It follows that we must have
$ {\displaystyle
\;\!
\mathbb{U}_{\!\;\!q}(t_0\:\!; t)
\leq
\lambda_{q}(t)
} $
in this case,
and the proof of Theorem 2.3 is complete.
\hfill $\Box$ \\
}
%
% --------------------------------------- End of proof
%                                        for Theorem 2.3
%
\mbox{} \vspace{-0.570cm} \\
\par
An important application of Theorem 2.3
is the following result. \\
\nl
%
% ------------------------------------------------------- %
%                                                         %
%                      Theorem 2.4                        %
%                                                         %
% ------------------------------------------------------- %
%
{\bf Theorem 2.4.}
\textit{%
Let
$\, p_{\mbox{}_{\!\;\!0}} \!\;\!\leq p < \infty $,
$\;\! 0 \leq t_0 \!\;\!< \mbox{\small $T$}_{\!\ast} $.
Then
} \\
\mbox{} \vspace{-0.750cm} \\
\begin{equation}
\tag{2.12}
\|\, u(\cdot,t) \,
\|_{\mbox{}_{\scriptstyle L^{\infty}(\mathbb{R})}}
\;\!\leq\:
\bigl(\;\! 2 \:\! p \;\!
\bigr)^{\mbox{}^{\scriptstyle \!\!\:\!\frac{1}{\scriptstyle p}}}
\!\cdot\;
\max\;\!
\biggl\{\,
\|\, u(\cdot,t_0) \,
\|_{\mbox{}_{\scriptstyle L^{\infty}(\mathbb{R})}}
\,\!;\;
\mathbb{B}(t_0\:\!; t)^{\mbox{}^{\scriptstyle \!\!\;\!\frac{1}{\scriptstyle p}}}
\!\:
\mathbb{U}_{\!\;\!p}(t_0\:\!; t)
\,\biggr\}
\end{equation}
\mbox{} \vspace{-0.200cm} \\
\textit{%
for any
$ \;\!t_0 \leq t < \mbox{\small $T$}_{\!\ast} $,
$\;\!$where
$ {\displaystyle
\;\!
\mathbb{B}(t_0\:\!; t)
} $,
$ {\displaystyle
\mathbb{U}_{\!\;\!p}(t_0\:\!; t)
\;\!
} $
are
given in $\;\!(2.9)$, $(2.10)$
above.
} \\
%
%
% -------------------------------------------------------
%
\mbox{} \vspace{+0.020cm} \\
{\small
{\bf Proof.}
Let $\;\!k \in \mathbb{Z} $, $ k \geq 2 $.
Applying (2.11) successively with
$ {\displaystyle
\;\!q \:\!=\:\!
2 \:\!p, \;\! 4\:\!p, ...\:\!, \;\! 2^{k}p
} $,
we obtain \\
\mbox{} \vspace{-0.200cm} \\
\mbox{} \hfill
$ {\displaystyle
\|\, u(\cdot,t) \,
\|_{\mbox{}_{\scriptstyle L^{2^{k}\!p}(\mathbb{R})}}
\!\:\!\leq\;
\max \,\biggl\{\,
\|\, \mbox{\boldmath $u$}(\cdot,t_0) \,
\|_{\mbox{}_{\scriptstyle L^{2^{k}p}(\mathbb{R})}}
\!\:\!;\:
K\!\:\!(k,\ell)^{\mbox{}^{\scriptstyle \!\:\!\frac{1}{\scriptstyle p}}}
\!\!\!\cdot\;\!
\mathbb{B}(t_0\:\!; t)^{\mbox{}^{\scriptstyle \!\:\!\frac{1}{\scriptstyle p}
\bigl(\,\!2^{\mbox{}^{\!-\ell}} \!\!\!\!\;\!-\; 2^{\mbox{}^{\!-k}}\,\!\bigr)}}
\hspace{-1.100cm} \cdot \hspace{0.770cm}
\|\, u(\cdot,t_0) \,
\|_{\mbox{}_{\scriptstyle L^{2^{\ell}\!p}(\mathbb{R})}}
\!\,\!,
} $ \\
\mbox{} \vspace{-0.500cm} \\
\mbox{} \hspace{9.50cm}
$ 1 \leq \ell \leq k - 1 \:\!; $ \\
\mbox{} \vspace{-0.450cm} \\
\mbox{} \hspace{+4.385cm}
$ {\displaystyle
K\!\:\!(k,0)^{\mbox{}^{\scriptstyle \!\:\!\frac{1}{\scriptstyle p}}}
\!\!\!\:\!\cdot\;\!
\mathbb{B}(t_0\:\!; t)^{\mbox{}^{\scriptstyle \!\:\!\frac{1}{\scriptstyle p}
\bigl(\:\!1\;\!-\: 2^{\mbox{}^{\!-k}}\,\!\bigr)}}
\hspace{-0.990cm} \cdot \hspace{0.700cm}
\mathbb{U}_{\!\;\!p}(t_0\:\!; t)
\,\biggr\}\,\!
} $,
\mbox{} \hfill (2.13$a$) \\
\mbox{} \vspace{-0.475cm} \\
where \\
\mbox{} \vspace{-0.475cm} \\
\mbox{} \hspace{+3.200cm}
$ {\displaystyle
K\!\:\!(k,\ell) \;\,=\,
\hspace{-0.250cm}
\prod_{\mbox{}\;\;j\,=\,\ell\;\!+\;\!1}^{k}
\!\!\!\!\;\!
\bigl(\;\! 2^{\mbox{}^{\scriptstyle \:\!j - 1}}
\!\:\!p \;C_{\mbox{}_{\!2}}^{\mbox{}^{\scriptstyle \;\!3}}
\;\!\bigr)^{\mbox{}^{\scriptstyle \!\,\! 2^{\mbox{}^{\!\;\!-\;\!j}} }}
\!\!\!\!\!,
\mbox{} \hspace{0.550cm}
0 \leq \ell \leq k - 1
} $.
\mbox{} \hfill (2.13$b$) \\

\mbox{} \vspace{-0.750cm} \\
Now,
for $\;\!1 \leq \ell \leq k - 1 $: \\
\mbox{} \vspace{-0.250cm} \\
\mbox{} \hspace{+2.250cm}
$ {\displaystyle
\mathbb{B}(t_0\:\!; t)^{\mbox{}^{\scriptstyle \!\:\!\frac{1}{\scriptstyle p}
\bigl(\,\!2^{\mbox{}^{\!-\ell}} \!\!\!\!\:\!-\; 2^{\mbox{}^{\!-k}}\,\!\bigr)}}
\hspace{-1.040cm} \cdot \hspace{0.770cm}
\|\, u(\cdot,t_0) \,
\|_{\mbox{}_{\scriptstyle L^{2^{\ell}\!p}(\mathbb{R})}}
} $ \\
\mbox{} \vspace{-0.250cm} \\
\mbox{} \hspace{+4.225cm}
$ {\displaystyle
\leq\;
\mathbb{B}(t_0\:\!; t)^{\mbox{}^{\scriptstyle \!\:\!\frac{1}{\scriptstyle p}
\bigl(\,\!2^{\mbox{}^{\!-\ell}} \!\!\!\!\:\!-\; 2^{\mbox{}^{\!-k}}\,\!\bigr)}}
\hspace{-1.040cm} \cdot \hspace{0.770cm}
\|\, u(\cdot,t_0) \,
\|_{\mbox{}_{\scriptstyle L^{p}(\mathbb{R})}}
  ^{\mbox{}^{\scriptstyle \!
    \frac{\scriptstyle \;\!
          2^{\mbox{}^{\!-\ell}} \!\!\!-\; 2^{\mbox{}^{\!-k}} }
          {\scriptstyle 1 \;\!-\: 2^{\mbox{}^{\!-k}} } }}
\hspace{-0.350cm} \cdot \;\,
\|\, u(\cdot,t_0) \,
\|_{\mbox{}_{\scriptstyle L^{2^{k}\!p}(\mathbb{R})}}
  ^{\mbox{}^{\scriptstyle \!
    \frac{\scriptstyle \;\!1 \;\!-\: 2^{\mbox{}^{\!-\ell}} }
         {\scriptstyle \;\!1 \;\!-\: 2^{\mbox{}^{\!-k}} } }}
} $ \\
\mbox{} \vspace{-0.050cm} \\
%
% \mbox{} \hspace{+4.550cm}
% $ {\displaystyle
% \leq\;
% \mathbb{B}(t_0\:\!; t)^{\mbox{}^{\scriptstyle \!\:\!\frac{1}{\scriptstyle p}
% \bigl(\:\!1 \;\!-\, 2^{\mbox{}^{\!-k}}\,\!\bigr)}}
% \hspace{-1.000cm} \cdot \hspace{0.670cm}
% \|\, u(\cdot,t_0) \,
% \|_{\mbox{}_{\scriptstyle L^{p}(\mathbb{R})}}
% \,+\;\,
% \|\, u(\cdot,t_0) \,
% \|_{\mbox{}_{\scriptstyle L^{2^{k}\!p}(\mathbb{R})}}
% } $ \\
% %
% \mbox{} \vspace{-0.100cm} \\
% %
\mbox{} \hspace{+4.250cm}
$ {\displaystyle
\leq\;
\max\,\biggl\{\,
\|\, u(\cdot,t_0) \,
\|_{\mbox{}_{\scriptstyle L^{2^{k}\!p}(\mathbb{R})}}
\!\:\!;\;\;\!
\mathbb{B}(t_0\:\!; t)^{\mbox{}^{\scriptstyle \!\:\!\frac{1}{\scriptstyle p}
\bigl(\:\!1 \;\!-\, 2^{\mbox{}^{\!-k}}\,\!\bigr)}}
\hspace{-1.000cm} \cdot \hspace{0.670cm}
\|\, u(\cdot,t_0) \,
\|_{\mbox{}_{\scriptstyle L^{p}(\mathbb{R})}}
\;\!\biggr\}
} $ \\
\mbox{} \vspace{+0.050cm} \\
by Young's inequality
(see e.g.$\;$\cite{Evans2002}, p.$\;$622);
in particular, we get, from (2.13), \\
\mbox{} \vspace{-0.150cm} \\
\mbox{} \hspace{+1.350cm}
% \mbox{} \hspace{+0.000cm}
$ {\displaystyle
\|\, u(\cdot,t) \,
\|_{\mbox{}_{\scriptstyle L^{2^{k}\!p}(\mathbb{R})}}
\!\;\!\leq\;
\bigl(\;\! 2 \:\! p \,
\bigr)^{\mbox{}^{\scriptstyle \!\!\:\!\frac{1}{\scriptstyle p}}}
\!\!\!\:\!\cdot\,
\max\,
\biggl\{\,
\|\, u(\cdot,t_0) \,
\|_{\mbox{}_{\scriptstyle L^{2^{k}\!p}(\mathbb{R})}}
\!\:\!;\;\,
\mathbb{B}(t_0\:\!; t)^{\mbox{}^{\scriptstyle \!\:\!\frac{1}{\scriptstyle p}
\bigl(\:\!1 \;\!-\, 2^{\mbox{}^{\!-k}}\,\!\bigr)}}
\hspace{-1.000cm} \cdot \hspace{0.670cm}
\mathbb{U}_{p}(t_0\:\!; t)
\,\biggr\}
} $, \\
\mbox{} \vspace{-0.000cm} \\
since
$ {\displaystyle
\:\!
K\!\;\!(k,\ell)
\leq
2 \:\! p \,
} $
for all
$ \;\! 0 \leq \ell \leq k - 1 $.
Letting $ k \rightarrow \infty $,
(2.12) is obtained.
\mbox{} \hfill $\Box$ \\
}
%
% --------------------------------------- End of proof
%                                        for Theorem 2.4
%
\mbox{} \vspace{-0.550cm} \\
\par
It follows from Theorems 2.1 and 2.4
that $ u(\cdot,t) $
is globally defined
($ \mbox{\small $T$}_{\!\ast} \!= \infty $).
Now,
from (2.12),
we immediately obtain,
letting $ \;\! t \rightarrow \infty $, \\
\mbox{} \vspace{-0.225cm} \\
\mbox{} \hspace{+0.315cm}
$ {\displaystyle
\limsup_{t\,\rightarrow\,\infty} \;\!
\|\, u(\cdot,t) \,
\|_{\mbox{}_{\scriptstyle L^{\infty}(\mathbb{R})}}
\leq\;\!
\bigl(\;\! 2 \;\! p \;\!
\bigr)^{\mbox{}^{\scriptstyle \!\!\:\!\frac{1}{\scriptstyle p}}}
\!\!\cdot\;
%
% \Bigr\{\,
% \|\, u(\cdot,t_0) \,
% \|_{\mbox{}_{\scriptstyle L^{\infty}(\mathbb{R})}}
% +\:
% \mathbb{B}(t_0)^{\mbox{}^{\scriptstyle \!\!\;\!\frac{1}{\scriptstyle p}}}
% \!\,
% \mathbb{U}_{\!\;\!p}(t_0)
% \;\!\Bigr\}
%
\max\,
\biggr\{\;\!
\|\, u(\cdot,t_0) \,
\|_{\mbox{}_{\scriptstyle L^{\infty}(\mathbb{R})}}
\!\:\!; \;\;\!
\mathbb{B}(t_0)^{\mbox{}^{\scriptstyle \!\!\;\!\frac{1}{\scriptstyle p}}}
\!\,
\mathbb{U}_{\!\;\!p}(t_0)
\;\!\biggr\}
} $
\hfill (2.14) \\
\mbox{} \vspace{+0.075cm} \\
for any
$ t_0 \geq 0 $,
where
$ \mathbb{B}(t_0) $,
$ \mathbb{U}_{\!\;\!p}(t_0) $
are given by \\
\mbox{} \vspace{-0.525cm} \\
\begin{equation}
\tag{2.15}
\mathbb{B}(t_0)
\;=\;\;\!
\sup\:\Bigl\{\;\!
% \mbox{\normalsize $ {\displaystyle \frac{\;\!B(t)}{\;\!\mu(t)} }$}
B(t)
\!\;\!:
\; t \geq t_0 \,\Bigr\},
\end{equation}
\mbox{} \vspace{-0.850cm} \\
\begin{equation}
\tag{2.16}
\mbox{} \;\;
\mathbb{U}_{p}(t_0)
\;=\;\;\!
\sup\:\Bigl\{\,
\|\, u(\cdot,t) \,
\|_{\mbox{}_{\scriptstyle L^{p}(\mathbb{R})}}
\!\!\:\!:
\; t \geq t_0 \,\Bigr\}.
\end{equation}
\mbox{} \vspace{-0.150cm} \\
Taking
% $ \:\!t_0^{\mbox{}} \!=\,\! t_0^{(n)} $
$ (\:\! t_0^{(n)} )_{n} \!\;\!$
such that
$ \:\!t_0^{(n)} \!\rightarrow \infty \,$
and
$ {\displaystyle
\,
\|\, u(\cdot,t_0^{(n)}) \,
\|_{\mbox{}_{\scriptstyle L^{\infty}(\mathbb{R})}}
\!\rightarrow\,
\liminf_{t\,\rightarrow\,\infty} \,
\|\, u(\cdot,t) \,
\|_{\mbox{}_{\scriptstyle L^{\infty}(\mathbb{R})}}
\!\:\!
} $, \linebreak
\mbox{} \vspace{-0.570cm} \\
and applying (2.14)
with $ \:\!t_0^{\mbox{}} \!=\,\! t_0^{(n)} $ for each $n$,
we then obtain,
letting $ n \rightarrow \infty $, \\
\mbox{} \vspace{-0.100cm} \\
\mbox{} \hspace{+0.300cm}
$ {\displaystyle
\limsup_{t\,\rightarrow\,\infty} \;\!
\|\, u(\cdot,t) \,
\|_{\mbox{}_{\scriptstyle L^{\infty}(\mathbb{R})}}
\!\;\!\leq\;\!
\bigl(\;\! 2 \:\! p \;\!
\bigr)^{\mbox{}^{\scriptstyle \!\!\:\!\frac{1}{\scriptstyle p}}}
\!\!\!\;\!\cdot\;\!
\max\,
\biggr\{\!\;\!
\liminf_{t\,\rightarrow\,\infty} \:\!
\|\, u(\cdot,t) \,
\|_{\mbox{}_{\scriptstyle L^{\infty}(\mathbb{R})}}
\!\;\!; \;\;\!
{\cal B}^{\mbox{}^{\scriptstyle \:\!\frac{1}{\scriptstyle p}}}
\!\!\!\,\!\cdot\;\!\!\;\!
{\cal U}_{p}
\;\!\biggr\},
} $
\hfill (2.17) \\
\mbox{} \vspace{+0.150cm} \\
where
$ \;\!{\cal B} $, $ {\cal U}_{p} $
are given by \\
\mbox{} \vspace{-0.700cm} \\
\begin{equation}
\tag{2.18}
{\cal B}
\;\!=\:
\limsup_{t\,\rightarrow\,\infty}
\,
% \mbox{\small $ {\displaystyle \frac{\small B(t)}{\small \mu(t)} }$}
B(t),
\qquad
{\cal U}_{p}
\;\!=\:
\limsup_{t\,\rightarrow\,\infty}
\;\!
\|\, u(\cdot,t) \,
\|_{\mbox{}_{\scriptstyle L^{p}(\mathbb{R})}}
\!\;\!.
\end{equation}
\nl
\mbox{} \vspace{-0.350cm} \\
%
%
% ************************************************************
% *                                                          *
% *         Section 3: The main large time estimates         *
% *                                                          *
% ************************************************************
%
\par
{\bf \S 3. Large time estimates} \\
\par
In this section,
we use the results obtained above
to derive two basic large time estimates
(given in Theorems 3.1 and 3.2 below)
for solutions $ u(\cdot,t) $
of problem (1.1$a$), (1.1$b$),
which represent important intermediate steps
that will ultimately lead to
the main result stated in Theorem 3.3. \\
\nl
%
% ------------------------------------------------------- %
%                                                         %
%                      Theorem 3.1                        %
%                                                         %
% ------------------------------------------------------- %
%
{\bf Theorem 3.1.}
\textit{%
Let
$ \;\!q \geq 2 \:\!p_{\mbox{}_{\!\;\!0}} \!\:\!$,
and
$ \;\!{\cal B} \!\;\!\geq 0 \:\!$
be as defined in $(2.18)$.
Then
} \\
\mbox{} \vspace{-0.725cm} \\
\begin{equation}
\tag{3.1}
\limsup_{t\,\rightarrow\,\infty} \;\!
\|\, u(\cdot,t) \,
\|_{\mbox{}_{\scriptstyle L^{q}(\mathbb{R})}}
\leq\,
\Bigl(\;\! \frac{\;\!q\;\!}{\mbox{\small $2$}}
\, C_{\mbox{}_{\!\;\!2}}^{\;\!3} \;\!
\Bigr)^{\scriptstyle \!\!
\frac{1}{\scriptstyle q} } \!\cdot\,
{\cal B}^{\mbox{}^{\scriptstyle \:\!
\frac{1}{\scriptstyle q} }} \!\cdot\;
\limsup_{t\,\rightarrow\,\infty} \;\!
\|\, u(\cdot,t) \,
\|_{\mbox{}_{\scriptstyle L^{q/2}(\mathbb{R})}}
\!\:\!,
% {\cal U}_{q/2}
\end{equation}
\mbox{} \vspace{-0.115cm} \\
\textit{%
where
$ {\displaystyle
\;\!C_{\mbox{}_{\!\;\!2}}
\!\:\!=\,\!
\bigl(\;\! 3 \;\!\sqrt{\:\!3\,} /\;\!(4 \:\!\pi) \:\!\bigr)^{\!1/3}
\!
} $
is the constant in the
Nash inequality $\;\!(2.7)$.
} \\
%
% -------------------------------------------------------
%
\mbox{} \vspace{-0.020cm} \\
{\small
{\bf Proof.}
We set $ \;\! p = q/2 \;\!$
and assume that
% $ {\cal B}_{\!\;\!\mu} $, $ {\cal U}_{p} $
$ {\cal U}_{p} $
is finite.
As in the proof of Theorem 2.2,
we take
$ {\displaystyle
\;\! v \in L^{\infty}(\mathbb{R} \times [\;\!0, \infty\:\![)
\;\!
} $
given by
$ {\displaystyle
v(x,t) = |\, u(x,t) \,|^{{\scriptstyle p}}
} $
if $ p > 1 $,
$ {\displaystyle
v(x,t) = u(x,t)
\;\!
} $
if $ p = 1 $.
It follows that \\
\mbox{} \vspace{-0.350cm} \\
\mbox{} \hspace{+4.500cm}
$ {\displaystyle
\|\, v(\cdot,t) \,
\|_{\mbox{}_{\scriptstyle L^{2}(\mathbb{R})}}^{2}
\!\;\!=\;
\|\, u(\cdot,t) \,
\|_{\mbox{}_{\scriptstyle L^{2p}(\mathbb{R})}}^{2p}
\!
} $, \\
\mbox{} \vspace{+0.100cm} \\
\mbox{} \hspace{+3.000cm}
$ {\displaystyle
\|\, v_{x}(\cdot,t) \,
\|_{\mbox{}_{\scriptstyle L^{2}(\mathbb{R})}}^{2}
\!\;\!=\;
p^{2} \!
\int_{\mathbb{R}} \!\;\!
|\, u(x,t) \,|^{\mbox{}^{\scriptstyle \;\!2\;\!
\mbox{\footnotesize $p$} \;\!-\;\! 2}}
\;\!
% \bigl|\, \frac{\partial \:\!u}{\partial \:\!x} \,\bigr|^{2}
|\, u_x(x,t) \,|^{\mbox{}^{\scriptstyle \;\!2}}
\, dx
} $. \\
\mbox{} \vspace{+0.100cm} \\
Therefore,
from (2.6),
we have,
for some null set
$ {\displaystyle
E_{\mbox{}_{2\:\!\mbox{\scriptsize $p$}}}
\!\subset [\;\!0, \infty \;\![
} $, \\
\mbox{} \vspace{-0.080cm} \\
\mbox{} \hspace{+1.000cm}
$ {\displaystyle
\frac{d}{d\:\!t} \,
\|\, v(\cdot,t) \,
\|_{\mbox{}_{\scriptstyle L^{2}(\mathbb{R})}}^{2}
\:\!+\;
4 \,\Bigl( 1 - \frac{1}{2\:\!p} \Bigr)
\,
\|\, v_{x}(\cdot,t) \,
\|_{\mbox{}_{\scriptstyle L^{2}(\mathbb{R})}}^{2}
% \:\!\leq
} $ \\
\mbox{} \vspace{+0.050cm} \\
\mbox{} \hfill
$ {\displaystyle
\leq\;
4 \,p\,
\Bigl( 1 - \frac{1}{2\:\!p} \Bigr)
\, B(t) \:
\|\, v(\cdot,t) \,
\|_{\mbox{}_{\scriptstyle L^{2}(\mathbb{R})}}
\:\!
\|\, v_{x}(\cdot,t) \,
\|_{\mbox{}_{\scriptstyle L^{2}(\mathbb{R})}}
% \!\;\!,
} $ \\
\mbox{} \vspace{-0.200cm} \\
for all
$ {\displaystyle
\;\! t \in [\;\!0, \infty\;\![ \,\setminus\;\!
E_{\mbox{}_{2\:\!\mbox{\scriptsize $p$}}}
} $,
and so,
by (2.7), \\
\mbox{} \vspace{-0.050cm} \\
\mbox{} \hspace{+1.000cm}
$ {\displaystyle
\frac{d}{d\:\!t} \,
\|\, v(\cdot,t) \,
\|_{\mbox{}_{\scriptstyle L^{2}(\mathbb{R})}}^{2}
\:\!+\;
4 \,\Bigl( 1 - \frac{1}{2\:\!p} \Bigr)
\,
\|\, v_{x}(\cdot,t) \,
\|_{\mbox{}_{\scriptstyle L^{2}(\mathbb{R})}}^{2}
} $ \\
\mbox{} \vspace{+0.050cm} \\
\mbox{} \hfill
$ {\displaystyle
\leq\;
4 \,p \:C_{\mbox{}_{\!\;\!2}} \:\!
\Bigl( 1 - \frac{1}{2\:\!p} \Bigr)
\;\! B(t) \:
\|\, v(\cdot,t) \,
\|_{\mbox{}_{\scriptstyle L^{1}(\mathbb{R})}}^{\;\!2/3}
\:\!
\|\, v_{x}(\cdot,t) \,
\|_{\mbox{}_{\scriptstyle L^{2}(\mathbb{R})}}^{\;\!4/3}
\!\:\!
} $. \\
\mbox{} \vspace{+0.050cm} \\
This gives,
by Young's inequality
(\cite{Evans2002}, p.$\;$622),
for all
$ {\displaystyle
\;\! t \in [\;\!0, \infty\;\![ \,\setminus\;\!
E_{\mbox{}_{2\:\!\mbox{\scriptsize $p$}}}
} $, \\
\mbox{} \vspace{-0.050cm} \\
\mbox{} \hspace{+1.000cm}
$ {\displaystyle
\frac{d}{d\:\!t} \,
\|\, v(\cdot,t) \,
\|_{\mbox{}_{\scriptstyle L^{2}(\mathbb{R})}}^{2}
\:\!+\;
\frac{\;\!4\;\!}{3} \;\!
\Bigl( 1 - \frac{1}{2\:\!p} \Bigr)
\,
\|\, v_{x}(\cdot,t) \,
\|_{\mbox{}_{\scriptstyle L^{2}(\mathbb{R})}}^{2}
\:\!\leq
} $ \\
\mbox{} \vspace{-0.850cm} \\
\mbox{} \hfill (3.2) \\
\mbox{} \vspace{-0.275cm} \\
\mbox{} \hfill
$ {\displaystyle
\leq\;
\frac{\;\!4\;\!}{3} \,
\Bigl( 1 - \frac{1}{2\:\!p} \Bigr)
\;\!
\bigl(\;\! p \:C_{\mbox{}_{\!\;\!2}} \:\!
\bigr)^{\!\;\!3}
\,
% \frac{\;\!B(t)^{3}}{\;\!\mu(t)^{2}} \;
B(t)^{3} \,
\|\, v(\cdot,t) \,
\|_{\mbox{}_{\scriptstyle L^{1}(\mathbb{R})}}^{\;\!2}
\!\:\!
} $. \\
\mbox{} \vspace{-0.150cm} \\
%
% 5555555555555555555555555555555555555555555
%
Setting \\
\mbox{} \vspace{-0.185cm} \\
\mbox{} \hspace{+1.400cm}
$ {\displaystyle
\lambda_{p} \;\!=\;
\limsup_{t\,\rightarrow\,\infty} \,
g(t),
\qquad
g(t) \,=\,
\bigl(\, p \:C_{\mbox{}_{\!\;\!2}}^{\;\!3}
\:\!\bigr)^{\mbox{}^{\scriptstyle \!\!\:\!1/2}}
% \Bigl(\;\! \frac{\;\!B(t)}{\;\mu(t)}
% \;\!\Bigr)^{\mbox{}^{\scriptstyle \!\!1/2}}
\!\;\!B(t)^{\mbox{}^{\scriptstyle \!\!1/2}} \,
\|\, v(\cdot,t) \,
\|_{\mbox{}_{\scriptstyle L^{1}(\mathbb{R})}}
\!\:\!
} $, \\
\mbox{} \vspace{-0.050cm} \\
we claim that \\
\mbox{} \vspace{-1.100cm} \\
\begin{equation}
\tag{3.3}
\limsup_{t\,\rightarrow\,\infty} \;\!
\|\, v(\cdot,t) \,
\|_{\mbox{}_{\scriptstyle L^{2}(\mathbb{R})}}
\leq\,
\lambda_{p}
\:\!.
\;\;\;
\end{equation}
\mbox{} \vspace{-0.200cm} \\
In fact, let us
argue by contradiction.
If (3.3) is false,
we can pick
$ {\displaystyle
\;\!0 < \eta \ll 1
\;\!
} $
and a sequence
$ (\;\!t_{j} \:\!)_{\mbox{}_{\scriptstyle \!\;\!j \,\geq\,0}} $,
$ t_{j} \rightarrow \infty $,
such that
$ {\displaystyle
\;\!
\|\, v(\cdot,t_{j}) \,
\|_{\mbox{}_{\scriptstyle L^{2}(\mathbb{R})}}
\!>\!\;\!
\lambda_{p} \!\;\!+ \eta
\,
} $
(for all $ j \geq 0 $)
and
$ {\displaystyle
\;\!
g(t) \leq
\lambda_{p} \!\;\!+ \eta/2
\;\!
} $
for all
$ \;\!t \geq t_0 $.
From (2.8$a$), Theorem 2.2,
it will then follow that \\
\mbox{} \vspace{-0.750cm} \\
\begin{equation}
\tag{3.4}
\|\, v(\cdot,t) \,
\|_{\mbox{}_{\scriptstyle L^{2}(\mathbb{R})}}
>\:\!
\lambda_{p} +\;\! \eta,
\qquad
\forall \;\, t \geq t_0 \:\!.
\end{equation}
\mbox{} \vspace{-0.275cm} \\
In fact,
suppose that
(3.4) were false,
so that
we had
$ {\displaystyle
\;\!
\|\, v(\cdot, \tilde{t}) \,
\|_{\mbox{}_{\scriptstyle L^{2}(\mathbb{R})}}
\!\!\:\!\leq\!\;\!
\lambda_{p} \!+ \eta
\;\!
} $
for some $ \;\!\tilde{t} > t_0 $.
Taking
$ j \gg 1 $
with
$ t_{j} \!> \tilde{t} $,
we could then find
$ \;\!\hat{t} \!\;\!\in [\,\tilde{t}, t_{j} \:\![ \:\! $
such that
$ {\displaystyle
\;\!
\|\, v(\cdot,t) \,
\|_{\mbox{}_{\scriptstyle L^{2}(\mathbb{R})}}
\!\!\;\!>\!\;\!
\lambda_{p} +\;\! \eta
\;\!
} $
for all
$ \;\!t \!\:\!\in \:]\,\hat{t}, t_{j} \:\!] $,
while
$ {\displaystyle
\;\!
\|\, v(\cdot,\hat{t}\:\!) \,
\|_{\mbox{}_{\scriptstyle L^{2}(\mathbb{R})}}
\!=
\lambda_{p} +\;\! \eta
} $,
and so
there would exist
$ {\displaystyle
t_{\ast} \!\in
[\,\hat{t}, t_{j} \:\!]
\:\!\setminus E_{\mbox{}_{2\;\!\mbox{\scriptsize $p$}}}
} $ \linebreak
with
$ {\displaystyle
\;\!
d/d\:\!t \,
\|\, v(\cdot,t) \,
\|_{\mbox{}_{\scriptstyle L^{2}(\mathbb{R})}}^{2}
\!
} $
positive at
$ \;\!t = t_{\ast} $.
By (2.8$a$),
we would have
$ {\displaystyle
\;\!
\|\, v(\cdot,t_{\ast}) \,
\|_{\mbox{}_{\scriptstyle L^{2}(\mathbb{R})}}
\!\leq
\lambda_{p}
} $,
but this would contradict
the fact that
$ {\displaystyle
\;\!
\|\, v(\cdot,t) \,
\|_{\mbox{}_{\scriptstyle L^{2}(\mathbb{R})}}
\!\geq
\lambda_{p} \!\,\!+ \eta
\,
} $
everywhere on
$ \;\![\, \hat{t}, t_{j} \,\!] $.
Thus,
we conclude that (3.4)
cannot be false,
as claimed.
$\!$We then obtain,
from (2.7), (3.2), (3.4), \linebreak
\mbox{} \vspace{-0.025cm} \\
\mbox{} \hspace{+1.750cm}
$ {\displaystyle
\|\, v(\cdot,t) \,
\|_{\mbox{}_{\scriptstyle L^{2}(\mathbb{R})}}^{\;\!6}
\;\!\leq\;
C_{\mbox{}_{\!2}}^{\;\!6} \,
\|\, v(\cdot,t) \,
\|_{\mbox{}_{\scriptstyle L^{1}(\mathbb{R})}}^{\;\!4}
\,
\|\, v_{x}(\cdot,t) \,
\|_{\mbox{}_{\scriptstyle L^{2}(\mathbb{R})}}^{\;\!2}
} $ \\
\mbox{} \vspace{+0.050cm} \\
\mbox{} \hspace{+4.165cm}
$ {\displaystyle
\leq\;
g(t)^{6}
\:+\;
\frac{2\:\!p}{\:\!2\:\!p - 1\:\!} \:
\|\, v(\cdot,t) \,
\|_{\mbox{}_{\scriptstyle L^{1}(\mathbb{R})}}^{\;\!4}
\!\;\!
\Bigl(\;\!- \; \frac{d}{d\:\!t} \,
\|\, v(\cdot,t) \,
\|_{\mbox{}_{\scriptstyle L^{2}(\mathbb{R})}}^{\;\!2}
\:\! \Bigr)
} $ \\
\mbox{} \vspace{+0.090cm} \\
for all
$ {\displaystyle
\, t \!\;\!\in [\;\!t_0, \infty\;\![ \,\setminus\;\!
E_{\mbox{}_{2\:\!\mbox{\scriptsize $p$}}}
} $.
Recalling that
$ {\displaystyle
\;\!
\|\, v(\cdot,t) \,
\|_{\mbox{}_{\scriptstyle L^{2}(\mathbb{R})}}
\!>
\lambda_{p} \!\;\!+ \eta
} $,
$ {\displaystyle
\:
g(t) \;\!\leq
\lambda_{p} \!\;\!+ \eta/2
\;\!
} $,
$ \; \forall \; t \geq t_0 $,
this gives \\
\mbox{} \vspace{-0.950cm} \\
\begin{equation}
\notag
- \; \frac{d}{d\:\!t} \,
\|\, v(\cdot,t) \,
\|_{\mbox{}_{\scriptstyle L^{2}(\mathbb{R})}}^{\;\!2}
\geq\:
K\!\:\!(\eta),
\qquad
\forall \;\,
t \in [\,t_0, \infty\;\![ \;\setminus\,
E_{\mbox{}_{2\:\!\mbox{\scriptsize $p$}}}
\end{equation}
\mbox{} \vspace{-0.200cm} \\
for some constant $ K\!\:\!(\eta) > 0\;\! $
independent of $\;\!t$,
which cannot be,
since this implies \\
\mbox{} \vspace{-0.650cm} \\
\begin{equation}
\notag
\|\, v(\cdot,t_0) \,
\|_{\mbox{}_{\scriptstyle L^{2}(\mathbb{R})}}^{\;\!2}
\geq\:
K\!\:\!(\eta) \cdot (\;\!t - t_0)
\qquad
\forall \;\,
t > t_0 \:\!.
\end{equation}
\mbox{} \vspace{-0.350cm} \\
This contradiction shows (3.3),
which is equivalent to (3.1),
and the proof is complete.
\mbox{} \hfill $\Box$ \\
}
%
% --------------------------------------- End of proof
%                                        for Theorem 3.1
%
\mbox{} \vspace{-0.850cm} \\
\par
Applying (3.1) successively
with
$ {\displaystyle
q =\:\!
2\:\!p, 4\:\!p, ... \:\!, 2^{k}p
} $,
we get \\
\mbox{} \vspace{-0.920cm} \\
\begin{equation}
\tag{3.5}
% {\cal U}_{\mbox{}_{\scriptstyle 2^{k}\!\:\!p}}
\limsup_{t\,\rightarrow\,\infty} \;\!
\|\, u(\cdot,t) \,\|_{\mbox{}_{\scriptstyle L^{2^{k}\!\:\!p}(\mathbb{R})}}
\!\;\!\leq\,
% K\!\:\!(k\:\!;\:\!p, m)^{\mbox{}^{\scriptstyle
% \!\;\!\frac{1}{\scriptstyle p}}}
\biggl[\;\;\!
\prod_{j\,=\,1}^{k} \;\!
\bigl(\;\! 2^{j-1} \:\!p\;\!
\,C_{\mbox{}_{\!2}}^{\;\!3} \:\!
\bigr)^{\mbox{}^{\scriptstyle \!\!\;\! 2^{\mbox{}^{\!-\;\!j}}}}
\biggr]^{\mbox{}^{\scriptstyle \!\frac{1}{\scriptstyle p}}}
\!\!\!\cdot\,
{\cal B}^{\mbox{}^{\scriptstyle
\;\!\frac{1}{\scriptstyle p}
\!\;\!\bigl( 1 \;\!-\: 2^{-k}\bigr)}}
\hspace{-0.950cm} \cdot \hspace{0.640cm}
{\cal U}_{p}
% \limsup_{t\,\rightarrow\,\infty} \,
% \|\, u(\cdot,t) \,\|_{\mbox{}_{\scriptstyle L^{p}(\mathbb{R})}}
\end{equation}
\mbox{} \vspace{-0.050cm} \\
for $ k \geq 1 $ arbitrary,
where
$ {\displaystyle
\;\!
{\cal U}_{p} \!\;\!=\;\!
\limsup_{t\,\rightarrow\,\infty} \;\!
\|\, u(\cdot,t) \,\|_{\mbox{}_{\scriptstyle L^{p}(\mathbb{R})}}
\!\;\!
} $.
Letting $ k \rightarrow \infty$,
this suggests \\
\mbox{} \vspace{-0.600cm} \\
\begin{equation}
\tag{3.6$a$}
% {\cal U}_{\mbox{}_{\scriptstyle 2^{k}\!\:\!p}}
\limsup_{t\,\rightarrow\,\infty} \;\!
\|\, u(\cdot,t) \,\|_{\mbox{}_{\scriptstyle L^{\infty}(\mathbb{R})}}
\leq\,
K\!\:\!(p)
\,\!\cdot\;\!
{\cal B}^{\mbox{}^{\scriptstyle
\frac{1}{\scriptstyle p}}}
\hspace{-0.150cm} \cdot \hspace{0.060cm}
% {\cal U}_{p},
\limsup_{t\,\rightarrow\,\infty} \;\!
\|\, u(\cdot,t) \,\|_{\mbox{}_{\scriptstyle L^{p}(\mathbb{R})}}
\!\;\!,
\end{equation}
\mbox{} \vspace{-0.100cm} \\
where \\
\mbox{} \vspace{-1.450cm} \\
\begin{equation}
\tag{3.6$b$}
K\!\:\!(p)
\;=\;
\biggl[\;\;\!
\prod_{j\,=\,1}^{\infty} \;\!
\bigl(\;\! 2^{j-1} \:\!p\;\!
\,C_{\mbox{}_{\!2}}^{\;\!3} \:\!
\bigr)^{\mbox{}^{\scriptstyle \!\!\;\! 2^{\mbox{}^{\!-\;\!j}}}}
\biggr]^{\mbox{}^{\scriptstyle \!\frac{1}{\scriptstyle p}}}
\!=\;
\Bigl(\;\! \frac{\mbox{\small $\;\!3 \;\! \sqrt{\,\!3\;\!} \;$}}
                {\mbox{\small $ 2 \:\! \pi $} } \: p \;\!
\Bigr)^{\mbox{}^{\scriptstyle \!\!
\frac{\scriptstyle 1}{\scriptstyle p} }}
\!,
\end{equation}
\mbox{} \vspace{-0.100cm} \\
%
% as given in (1.6),
cf.$\;$(1.6) above,
as long as
the limit processes
$ k \rightarrow \infty$, $ t \rightarrow \infty $
can be interchanged.
That this is indeed the case
is a consequence of (2.17) and
the following result. \\
\nl
%
% ------------------------------------------------------- %
%                                                         %
%                      Theorem 3.2                        %
%                                                         %
% ------------------------------------------------------- %
%
{\bf Theorem 3.2.}
\textit{%
Let
$ \;\!p \geq p_{\mbox{}_{\!\;\!0}} \!\:\!$.
Then
} \\
\mbox{} \vspace{-0.850cm} \\
\begin{equation}
\tag{3.7}
\liminf_{t\,\rightarrow\,\infty} \;\!
\|\, u(\cdot,t) \,
\|_{\mbox{}_{\scriptstyle L^{\infty}(\mathbb{R})}}
\leq\,
\bigl(\, p \,
C_{\mbox{}_{\!\;\!2}} \:\! C_{\mbox{}_{\!\!\;\!\infty}}
\bigr)^{\mbox{}^{\scriptstyle \!\!\;\!
\frac{1}{\scriptstyle p} }} \!\!\!\;\!\cdot\,
{\cal B}^{\mbox{}^{\scriptstyle \:\!
\frac{1}{\scriptstyle p} }} \!\cdot\;
\limsup_{t\,\rightarrow\,\infty} \;\!
\|\, u(\cdot,t) \,
\|_{\mbox{}_{\scriptstyle L^{p}(\mathbb{R})}}
\!\;\!,
% {\cal U}_{p}
\end{equation}
\mbox{} \vspace{-0.220cm} \\
\textit{%
where
$ {\displaystyle
\;\!C_{\mbox{}_{\!\;\!2}}, \, C_{\mbox{}_{\!\infty}}
\!\:\!
} $
are the constants given in $\;\!(2.7)$, $(1.9)$.
} \\
%
% -------------------------------------------------------
%
\mbox{} \vspace{-0.040cm} \\
{\small
{\bf Proof.}
Again,
assuming
$ {\displaystyle
\;\!
{\cal U}_{p} \!\;\!
} $
finite
(otherwise, (3.7) is obvious, cf.$\;$footnote 1),
we introduce,
as in the previous proof,
$ {\displaystyle
\;\!
v \in L^{\infty}(\mathbb{R} \times [\;\!0, \infty\:\![)
\;\!
} $
given by
$ {\displaystyle
v(x,t) = |\, u(x,t) \,|^{{\scriptstyle p}}
} $
if $ p > 1 $,
and
$ {\displaystyle
v(x,t) = u(x,t)
\;\!
} $
if $ p = 1 $.
Thus,
(3.2) is valid,
and setting
$ \lambda_{p} \!\in \mathbb{R} $,
$ \;\!g \in L^{\infty}([\;\!0, \infty\;\![\:\!) \;\!$
by \\
\mbox{} \vspace{-0.675cm} \\
\begin{equation}
\notag
\lambda_{p} \;\!=\;
\limsup_{t\,\rightarrow\,\infty} \,
g(t),
\qquad
g(t) \,=\,
p \:C_{\mbox{}_{\!\;\!2}} \,
% \frac{\;\!B(t)}{\;\mu(t)} \,
B(t) \,
\|\, \mbox{\boldmath $v$}(\cdot,t) \,
\|_{\mbox{}_{\scriptstyle L^{1}(\mathbb{R})}}
\!\:\!,
\end{equation}
\mbox{} \vspace{-0.150cm} \\
we have that (3.7) is obtained
if we show that \\
\mbox{} \vspace{-0.700cm} \\
\begin{equation}
\tag{3.8}
\liminf_{t\,\rightarrow\,\infty} \;\!
\|\, v(\cdot,t) \,
\|_{\mbox{}_{\scriptstyle L^{\infty}(\mathbb{R})}}
\;\!\leq\;
C_{\mbox{}_{\!\infty}} \!\!\cdot \lambda_{p}
\:\!.
\;\;\;
\end{equation}
\mbox{} \vspace{-0.500cm} \\
We argue by contradiction
and assume that (3.8) is false.
Taking then
$ {\displaystyle
\;\!0 < \eta \ll 1,
\; t_0 \gg 1
} $
so that
$ {\displaystyle
\;\!
\|\, v(\cdot,t) \,
\|_{\mbox{}_{\scriptstyle L^{\infty}(\mathbb{R})}}
\!\geq\;\!
C_{\mbox{}_{\!\infty}} \!\!\cdot
(\lambda_{p} \!\;\!+ \eta \:\!)
\,
} $
and
$ {\displaystyle
\:
g(t) \;\!\leq
\lambda_{p} \!\;\!+ \eta/2
\;\!
} $
hold
for all
$ t \geq t_0 $,
we get, \linebreak
by (1.9), (3.2), \\
\mbox{} \vspace{-0.175cm} \\
\mbox{} \hspace{+0.750cm}
$ {\displaystyle
\|\, v(\cdot,t) \,
\|_{\mbox{}_{\scriptstyle L^{\infty}(\mathbb{R})}}^{\;\!3}
\;\!\leq\;
C_{\mbox{}_{\!\infty}}^{\;\!3} \,
\|\, v(\cdot,t) \,
\|_{\mbox{}_{\scriptstyle L^{1}(\mathbb{R})}}
\,
\|\, v_{x}(\cdot,t) \,
\|_{\mbox{}_{\scriptstyle L^{2}(\mathbb{R})}}^{\;\!2}
} $ \\
\mbox{} \vspace{+0.050cm} \\
\mbox{} \hfill
$ {\displaystyle
\leq\;
C_{\mbox{}_{\!\!\;\!\infty}}^{\,3} \: g(t)^{3}
\:+\;
C_{\mbox{}_{\!\!\;\!\infty}}^{\,3} \;
\frac{2\:\!p}{\:\!2\:\!p - 1\:\!} \;
\|\, v(\cdot,t) \,
\|_{\mbox{}_{\scriptstyle L^{1}(\mathbb{R})}}
\:\!
\Bigl(\;\!- \; \frac{d}{d\:\!t} \,
\|\, v(\cdot,t) \,
\|_{\mbox{}_{\scriptstyle L^{2}(\mathbb{R})}}^{\;\!2}
\:\! \Bigr)
} $ \\
\mbox{} \vspace{+0.090cm} \\
for all
$ {\displaystyle
\;\! t \in [\;\!t_0, \infty\;\![ \,\setminus\;\!
E_{\mbox{}_{2\:\!\mbox{\scriptsize $p$}}}
} $.
Since
$ {\displaystyle
\;\!
\|\, v(\cdot,t) \,
\|_{\mbox{}_{\scriptstyle L^{\infty}(\mathbb{R})}}
\!\geq\;\!
C_{\mbox{}_{\!\infty}} \!\!\cdot
(\lambda_{p} \!\;\!+ \eta \:\!)
} $,
$ {\displaystyle
\:
g(t) \;\!\leq
\lambda_{p} \!\;\!+ \eta/2
\;\!
} $,
$\:\!$this gives \\
\mbox{} \vspace{-0.550cm} \\
\begin{equation}
\notag
- \; \frac{d}{d\:\!t} \,
\|\, v(\cdot,t) \,
\|_{\mbox{}_{\scriptstyle L^{2}(\mathbb{R})}}^{\;\!2}
\;\!\geq\;
K\!\:\!(\eta),
\qquad
\forall \;\,
t \in [\,t_0, \infty\;\![ \:\setminus\,
E_{\mbox{}_{2\:\!\mbox{\scriptsize $p$}}}
\end{equation}
\mbox{} \vspace{-0.200cm} \\
for some constant $ K\!\;\!(\eta) > 0\;\! $
independent of $\,t$.
As before,
this implies
that
$ {\displaystyle
\,
\|\, v(\cdot,t_0) \,
\|_{\mbox{}_{\scriptstyle L^{2}(\mathbb{R})}}^{\;\!2}
\!\geq
} $
$ {\displaystyle
K\!\:\!(\eta) \cdot (\:\! t - t_0)
\;\!
} $
for all $ \;\! t \geq t_0 $,
which is impossible because
$ {\displaystyle
\|\, v(\cdot,t_0) \,
\|_{\mbox{}_{\scriptstyle L^{2}(\mathbb{R})}}
\!
} $
is finite.
This contradiction establishes
(3.8) above,
completing
the proof of Theorem 3.2.
\mbox{} \hfill $\Box$ \\
}
%
% --------------------------------------- End of proof
%                                        for Theorem 3.2
\mbox{} \vspace{-0.700cm} \\
\par
We are finally in good position
to derive (1.6), (3.6).
$\!$Combining (2.17) and (3.7) above,
we obtain \\
\mbox{} \vspace{-1.050cm} \\
\begin{equation}
\tag{3.9}
\limsup_{t\,\rightarrow\,\infty} \;\!
\|\, u(\cdot,t) \,
\|_{\mbox{}_{\scriptstyle L^{\infty}(\mathbb{R})}}
\leq\,
\bigl(\, 2 \;\! p^{2} \:\!
\bigr)^{\mbox{}^{\scriptstyle \!\!\;\!
\frac{1}{\scriptstyle p} }} \!\!\!\,\cdot\,
{\cal B}^{\mbox{}^{\scriptstyle \:\!
\frac{1}{\scriptstyle p} }} \!\cdot\;
{\cal U}_{p}
% \limsup_{t\,\rightarrow\,\infty} \;\!
% \|\, u(\cdot,t) \,
% \|_{\mbox{}_{\scriptstyle L^{p}(\mathbb{R})}}
\end{equation}
\mbox{} \vspace{-0.175cm} \\
for each
$ p \!\;\!\geq p_{\mbox{}_{\!\;\!0}} $,
so that we have,
in particular, \\
\mbox{} \vspace{-0.890cm} \\
\begin{equation}
\tag{3.10}
\limsup_{t\,\rightarrow\,\infty} \;\!
\|\, u(\cdot,t) \,
\|_{\mbox{}_{\scriptstyle L^{\infty}(\mathbb{R})}}
\leq\,
\bigl(\;\! 2^{2\:\! k \,+\, 1} \;\!p^{2}
\bigr)^{\mbox{}^{\scriptstyle \!\!\;\!
\frac{1}{\scriptstyle 2^{k}\!\:\!p} }}
\hspace{-0.380cm} \cdot\;
{\cal B}^{\mbox{}^{\scriptstyle \:\!
\frac{1}{\scriptstyle 2^{k}\!\:\!p} }} \!\cdot\;
{\cal U}_{2^{k}\!\:\!p}
% \limsup_{t\,\rightarrow\,\infty} \;\!
% \|\, \mbox{\boldmath $u$}(\cdot,t) \,
% \|_{\mbox{}_{\scriptstyle L^{p}(\mathbb{R})}}
\end{equation}
\mbox{} \vspace{-0.125cm} \\
for each $ k \geq 0 $.
By (3.5),
we then get \\
\mbox{} \vspace{-0.975cm} \\
\begin{equation}
\tag{3.11}
\limsup_{t\,\rightarrow\,\infty} \;\!
\|\, u(\cdot,t) \,
\|_{\mbox{}_{\scriptstyle L^{\infty}(\mathbb{R})}}
\leq\,
\biggl\{\;\!
\bigl(\;\! 2^{2\:\! k \,+\, 1} \;\!p^{2}
\bigr)^{\mbox{}^{\scriptstyle \!\!\;\!2^{\mbox{}^{\!-\;\!k}} }}
\hspace{-0.400cm} \:\!\cdot\;
\prod_{j\,=\,1}^{k} \;\!
\bigl(\;\! 2^{j-1} \:\!p\;\!
\,C_{\mbox{}_{\!2}}^{\;\!3} \:\!
\bigr)^{\mbox{}^{\scriptstyle \!\!\;\! 2^{\mbox{}^{\!-\;\!j}}}}
\biggr\}^{\mbox{}^{\scriptstyle \!\!\frac{1}{\scriptstyle p} }}
\hspace{-0.150cm} \cdot\:
{\cal B}^{\mbox{}^{\scriptstyle \:\!
\frac{1}{\scriptstyle p} }} \!\cdot\;
{\cal U}_{p}
% \limsup_{t\,\rightarrow\,\infty} \;\!
% \|\, \mbox{\boldmath $u$}(\cdot,t) \,
% \|_{\mbox{}_{\scriptstyle L^{p}(\mathbb{R})}}
\end{equation}
\mbox{} \vspace{-0.075cm} \\
for all $\;\! k $.
Letting $ \;\!k \rightarrow \infty $,
% (1.6) is obtained.
Theorem 3.3 is obtained,
and our argument is complete. \\
\nl
%
% ------------------------------------------------------- %
%                                                         %
%                      Theorem 3.3                        %
%                                                         %
% ------------------------------------------------------- %
%
{\bf Theorem 3.3.}
\textit{%
Let
$ \;\!p \geq p_{\mbox{}_{\!\;\!0}} \!\;\!$.
Assuming
$ {\displaystyle
\;\!b \in L^{\infty}(\mathbb{R} \times [\;\!0, \infty\;\![\;\!)
} $,
then $(1.6)$, $(3.6)$ hold.
} \\
\mbox{} \vspace{-0.550cm} \\
\par
It is worth noticing that
the corresponding estimate
for the $n$-dimensional problem (1.8),
namely, \\
\mbox{} \vspace{-0.900cm} \\
\begin{equation}
\tag{3.12}
\limsup_{t\,\rightarrow\,\infty}\;\!
\|\, u(\cdot,t) \,
\|_{\mbox{}_{\scriptstyle L^{\infty}(\mathbb{R}^{n})}}
\leq\,
K\!\:\!(n,p) \;\!\cdot\;\!
{\cal B}^{\mbox{}^{\scriptstyle
\;\!\frac{\scriptstyle n}{\scriptstyle p} }}
\!\!\cdot\:
% {\cal U}_{p}
\limsup_{t\,\rightarrow\,\infty}\;\!
\|\, u(\cdot,t) \,
\|_{\mbox{}_{\scriptstyle L^{p}(\mathbb{R}^{n})}}
\!\:\!,
\end{equation}
\vspace{-0.200cm} \\
where $ {\cal B} \geq 0 $ is similarly defined,
can be also derived in arbitrary dimension $ n > 1 $.

\mbox{} \vspace{-0.550cm} \\
%
%
% ************************************************************
% *                                                          *
% *                  Section 4: Open problems                *
% *                                                          *
% ************************************************************
%
\par
{\bf \S 4. Concluding remarks} \\
\par
We close our discussion of the problem (1.1$a$), (1.1$b$),
given
$ {\displaystyle
\;\!b \in L^{\infty}(\mathbb{R} \times [\;\!0, \infty\;\![\;\!)
} $,
$ 1 \leq p_{\mbox{}_{\!\;\!0}} \!< \infty $,
indicating a few questions
which were not answered
by our analysis: \\
\mbox{} \vspace{-0.100cm} \\
({\em a\/})
\begin{minipage}[t]{14.30cm}
characterize all
$ {\displaystyle
\;\!b \in L^{\infty}(\mathbb{R} \times [\;\!0, \infty\;\![\;\!)
} $
for which it is true that
$ {\displaystyle
\,
\|\, u(\cdot,t) \,\|_{\mbox{}_{\scriptstyle L^{\infty}(\mathbb{R})}}
\!\!\rightarrow 0
} $
(as $ t \rightarrow \infty $)
for every solution $ u(\cdot,t) $ of problem (1.1);
\end{minipage}
\nl
\mbox{} \vspace{-0.350cm} \\
({\em b\/})
\begin{minipage}[t]{14.30cm}
same question as ({\em a\/}) above,
but requiring only that
$ {\displaystyle
\;\!
\limsup
\|\, u(\cdot,t) \,\|_{\mbox{}_{\scriptstyle L^{\infty}(\mathbb{R})}}
\!\!\;\!< \infty
} $
(as $ t \rightarrow \infty $)
for every solution $ u(\cdot,t) $ of problem (1.1),
in case $ \;\!p_{\mbox{}_{\!\;\!0}} \!> 1 $;\footnotemark
\end{minipage}
\footnotetext{%
$\:\!$For $ p_{\mbox{}_{\!\;\!0}} \!= 1 $,
any
$ {\displaystyle
\;\!b \in L^{\infty}(\mathbb{R} \times [\;\!0, \infty\;\![\;\!)
} $
satisfies property ({\em b\/}), cf.$\;$(1.7) in Section 1. \\
\mbox{} \vspace{-0.900cm} \\
}
\nl
\mbox{} \vspace{-0.350cm} \\
({\em c\/})
\begin{minipage}[t]{14.30cm}
given $ p_{\mbox{}_{\!\;\!0}} > 1 $,
characterize all
$ {\displaystyle
\;\!b \in L^{\infty}(\mathbb{R} \times [\;\!0, \infty\;\![\;\!)
} $
such that
$ {\displaystyle
\,
\|\, u(\cdot,t) \,\|_{\mbox{}_{\scriptstyle L^{p_{\mbox{}_{\!\;\!0}}}(\mathbb{R})}}
\!\!\rightarrow 0
} $
(as $ t \rightarrow \infty $)
for every solution $ u(\cdot,t) $ of problem (1.1);
\end{minipage}
\nl
\mbox{} \vspace{-0.300cm} \\
({\em d\/})
\begin{minipage}[t]{14.30cm}
same question as ({\em c\/}) above,
but requiring only that
$ {\displaystyle
\;\!
\limsup
\|\, u(\cdot,t) \,\|_{\mbox{}_{\scriptstyle L^{p_{\mbox{}_{\!\;\!0}}}(\mathbb{R})}}
\!\!\;\!< \infty
} $
(as $ t \rightarrow \infty $)
for every solution $ u(\cdot,t) $ of problem (1.1);
\end{minipage}
\nl
\mbox{} \vspace{-0.350cm} \\
({\em e\/})
\begin{minipage}[t]{14.30cm}
for $ p_{\mbox{}_{\!\;\!0}} = 1 $,
characterize all
$ {\displaystyle
\;\!b \in L^{\infty}(\mathbb{R} \times [\;\!0, \infty\;\![\;\!)
} $
such that
$ {\displaystyle
\,
\|\, u(\cdot,t) \,\|_{\mbox{}_{\scriptstyle L^{1}(\mathbb{R})}}
\!\!\rightarrow |\,m\,|
} $
(as $ \,\!t \rightarrow \infty $)
for every solution $ u(\cdot,t) $,
where $ \,\!m = \!\;\!\int_\mathbb{R} \!u_0(x)\;\!dx\;\! $ is the solution mass;
\end{minipage}
\nl
\mbox{} \vspace{-0.350cm} \\
({\em f\/})
\begin{minipage}[t]{14.30cm}
for $ p_{\mbox{}_{\!\;\!0}} = 1 $,
and
$ {\displaystyle
\;\!b \in L^{\infty}(\mathbb{R} \times [\;\!0, \infty\;\![\;\!)
} $
not satisfying property ({\em e\/}) above,
what are the values of
$ {\displaystyle
\lim_{t\,\rightarrow\,\infty} \;\!
\|\, u(\cdot,t) \,\|_{\mbox{}_{\scriptstyle L^{1}(\mathbb{R})}}
} $
in case of initial states that change sign?
\end{minipage}
\nl
\mbox{} \vspace{-0.300cm} \\
These questions can be similarly posed
for solutions $ u(\cdot,t) $
of autonomous problems \\
\mbox{} \vspace{-0.650cm} \\
\begin{equation}
\tag{4.1}
u_t \;\!+\, (\;\!b(x) \;\!u \;\!)_{x}
\;\!=\;
u_{xx},
\qquad
u(\cdot,0) \in
L^{p_{\mbox{}_{\!\;\!0}}}(\mathbb{R} \cap L^{\infty}(\mathbb{R})
\end{equation}
\mbox{} \vspace{-0.260cm} \\
where $ b \in L^{\infty}(\mathbb{R}) $
does not depend on the time variable.
For (4.1), question ({\em e\/})
has been answered in
\cite{Rudnicki1993}.
(See also \cite{BrzezniakSzafirski1991}).
Another interesting question
is the following: \\
\mbox{} \vspace{-0.125cm} \\
({\em g\/})
\begin{minipage}[t]{14.30cm}
when (4.1) admits no stationary solutions
other than the trivial solution $ u = 0 $, \linebreak
is it true that
$ {\displaystyle
\lim_{t\,\rightarrow\,\infty}
\|\, u(\cdot,t) \,\|_{\mbox{}_{\scriptstyle L^{\infty}(\mathbb{R})}}
\! =\;\! 0
\,
} $
for every solution $ u(\cdot,t) $?
\end{minipage}
\nl
\mbox{} \vspace{-0.375cm} \\
Moreover,
for solutions
$ u(\cdot,t) $
of (1.1) or (4.1)
with
$ {\displaystyle
\|\, u(\cdot,t) \,\|_{\mbox{}_{\scriptstyle L^{\infty}(\mathbb{R})}}
\!\! \rightarrow 0
\;\!
} $
as $ \;\!t \rightarrow \infty $,
there is the question
of determining the proper decay rate.\footnote{%
$\:\!$In case we have $ b_x \geq 0 $ for all $x$, $t$,
the answer is given in (1.2) above.
}
As suggested by Fig.$\,$1,
solution decay may sometimes happen
at remarkably slow rates.  \\
%
% ------------------------------------------------------------
%
\nl
\nl
%
%
% ************************************************************
% *                                                          *
% *                  Acknowledgments                         *
% *                                                          *
% ************************************************************
%
\mbox{} \vspace{-2.500cm} \\
{\bf Acknowledgements.}
The authors would like to thank {\tt CNPq}
(Conselho Nacional de Desenvolvimento Cient\'\i fico
 e Tecnol\'ogico, Brazil)
for their financial support. \\
\mbox{} \vspace{-1.000cm} \\
%
%
% ************************************************************
% *                                                          *
% *                      References                          *
% *                                                          *
% ************************************************************
%

%
% -----------------------------------------------------------
%

\nl
\mbox{} \vspace{-0.300cm} \\
%

%
% -----------------------------------------------------------
%
\nl
\begin{minipage}[t]{10.00cm}
\textsc{Jos\'e Afonso Barrionuevo} \\
Departamento de Matem\'atica Pura e Aplicada \\
Universidade Federal do Rio Grande do Sul \\
Porto Alegre, RS 91509-900, Brazil \\
E-mail: {\sf josea@mat.ufrgs.br}
\end{minipage}
\nl
\nl
\mbox{} \vspace{-0.150cm} \\
%
%\mbox{} \hspace{3.50cm}
%
\begin{minipage}[t]{10.00cm}
\textsc{Lucas da Silva Oliveira} \\
Departamento de Matem\'atica Pura e Aplicada \\
Universidade Federal do Rio Grande do Sul \\
Porto Alegre, RS 91509-900, Brazil \\
E-mail: {\sf lucas.oliveira@ufrgs.br}
\end{minipage}
\nl
\nl
\mbox{} \vspace{-0.150cm} \\
\begin{minipage}[t]{10.00cm}
\textsc{Paulo Ricardo Zingano} \\
Departamento de Matem\'atica Pura e Aplicada \\
Universidade Federal do Rio Grande do Sul \\
Porto Alegre, RS 91509-900, Brazil \\
E-mail: {\sf paulo.zingano@ufrgs.br}
\end{minipage}

%
% ----------------------------------------------------------------
%

\end{document}